\documentclass{elsarticle}

\usepackage{hyperref}

\journal{Computers \& Mathematics with Applications}

\usepackage{color}
\usepackage{graphicx}

\usepackage{lmodern}
\usepackage{mathrsfs}
\usepackage{multirow}
\usepackage{amsmath,amstext,amssymb}
\usepackage{amsthm}
\usepackage{tabularx,longtable,multirow,diagbox}
\usepackage{booktabs}

\usepackage{bm}

\usepackage[ruled,vlined,linesnumbered]{algorithm2e}
\SetKwInput{KwIn}{Input}
\SetKwInput{KwOut}{Output}
\SetKw{KwRet}{Return}
\DontPrintSemicolon

\numberwithin{equation}{section}

\theoremstyle{plain}
\newtheorem{theorem}{Theorem}[section]

\newtheorem{lemma}[theorem]{Lemma}

\theoremstyle{remark}

\SetAlgoNoLine                       
\SetNlSty{}{}{}                      
\SetAlgoNlRelativeSize{-1}          
\DontPrintSemicolon                 
\SetInd{1em}{1em}                  

\usepackage{booktabs}
\usepackage{subcaption}
\usepackage{enumitem}

\usepackage{xcolor}

\newcommand{\rd}{\,\mathrm{d}}
\newcommand{\dx}{\,\mathrm{d}x}
\newcommand{\dy}{\,\mathrm{d}y}
\newcommand{\ds}{\,\mathrm{d}s}

\usepackage{placeins}
\usepackage{needspace}

\begin{document}
	
\begin{frontmatter}
		
\title{A Two-stage Adaptive Lifting PINN Framework for Solving Viscous Approximations to Hyperbolic Conservation Laws}

\author[nju]{Yameng Zhu} 
\ead{YMZhu@smail.nju.edu.cn}
\author[nju]{Weibing Deng\corref{mycorrespondingauthor}}
\ead{wbdeng@nju.edu.cn}
\author[nju]{Ran Bi} 
\ead{ranbi@smail.nju.edu.cn}
\cortext[mycorrespondingauthor]{Corresponding author}
\address[nju]{School of Mathematics,
	Nanjing University, Nanjing 210093, People's Republic of China}


\begin{abstract}
Training physics-informed neural networks (PINNs) for hyperbolic conservation laws near the inviscid limit presents considerable difficulties: strong-form residuals are not classically defined at shock discontinuities, while small-viscosity regularization introduces narrow internal layers that exacerbate spectral bias. While lifting-based strategies can alleviate these representation difficulties by embedding the problem into a higher-dimensional space, a key practical challenge is the automatic construction of auxiliary variables that encode shock locations and interface geometry. To overcome this limitation, this paper proposes a novel Two-stage Adaptive Lifting PINN--a fully autonomous framework. The key idea is to augment the physical coordinates by introducing a learned auxiliary field generated through $r$-adaptive coordinate transformations. Crucially, this geometric manifold is automatically extracted from a computationally inexpensive, high-viscosity coarse solution in the first stage. This pre-conditioning strategy effectively guides the network to capture steep gradients at the target small viscosity without requiring any prior knowledge of the shock dynamics. The coordinate construction and lifted training can be repeated within a viscosity continuation scheme that progressively reduces the viscosity. Theoretically, we first derive a scalar a posteriori $L^2$ error estimate to quantify the viscosity dependence of residual-based error control. Secondly, we provide a statistical interpretation establishing a variance-reduction criterion for coordinate-induced importance sampling. Finally, we perform a neural tangent kernel analysis of the gradient-flow dynamics, explaining how adaptive augmentation enriches the tangent-feature representation and accelerates residual decay. Our numerical experiments demonstrate faster loss decay and improved small-viscosity accuracy for scalar equations and the Euler system. Ablation results show that lifting enables further accuracy gains as viscosity decreases, beyond those obtained by viscosity continuation and adaptive sampling alone.
\end{abstract}


\end{frontmatter}


\section{Introduction}
\label{sec:introduction}

Many physical and engineering systems are governed by partial differential equations (PDEs) that can develop sharp gradients or discontinuities over time. A prominent class of such equations is hyperbolic conservation laws~\cite{dafermos2005hyperbolic}, which include the Euler equations in fluid dynamics, the inviscid Burgers equation, and related models. Even with smooth initial and boundary conditions, their solutions may evolve into discontinuous structures such as shock waves and contact discontinuities, presenting substantial challenges for numerical methods.

Traditional numerical approaches--such as finite difference, finite volume~\cite{leveque2002finite}, and finite element methods~\cite{johnson1983finite}---typically discretize the domain using a mesh to approximate spatial derivatives. However, accurately capturing sharp gradients or discontinuities often demands extremely fine spatial resolution, which incurs high computational cost. To mitigate this, high-resolution schemes have been developed, including flux-corrected transport~\cite{boris1973flux}, total variation diminishing~\cite{harten1997high}, and weighted essentially non-oscillatory methods~\cite{jiang1996efficient}, which effectively suppress spurious oscillations near discontinuities. Additional strategies include adaptive mesh refinement~\cite{berger1989local}, which dynamically enhances local resolution, and discontinuous Galerkin methods~\cite{cockburn2001runge}, offering high-order accuracy and robust shock handling on unstructured grids.

Despite their successes, resolving narrow transition layers remains demanding, and adaptive resolution plays a central role in controlling the associated cost. These challenges also motivate learning-based approaches that adapt their representations to local solution features.

The rapid progress in computing and machine learning has spurred the development of neural-network (NN) based methods for PDE solution and flow prediction~\cite{pang2019neural, lye2020deep, magiera2020constraint}. Among these, PINNs~\cite{raissi2019pinn} have garnered significant interest, especially in fluid mechanics~\cite{cai2021pinns_review}. By incorporating governing equations and boundary conditions directly into the loss function through automatic differentiation, PINNs allow residual evaluation at freely chosen collocation points without a mesh-based differentiation stencil. Furthermore, PINNs provide a unified framework that accommodates both forward and inverse problems.

However, directly applying vanilla PINNs to inviscid hyperbolic conservation laws faces a fundamental difficulty: at shocks and contact discontinuities, the solution ceases to be classically differentiable, so the pointwise strong-form residual is not defined~\cite{liu2024discontinuity}. To address this difficulty, weak PINNs impose entropy conditions through weak residuals~\cite{deryck2024wpinns}, while localized linearizations modify the residual formulation near shocks in Riemann problems~\cite{liu2024locally}. Domain-decomposition methods such as Extended PINNs (XPINNs)~\cite{jagtap2020xpinns} use separate networks on space--time subdomains with interface coupling.

A common mitigation strategy is to introduce a small viscous regularization term, thereby transforming the system into a parabolic form~\cite{coutinho2023physics}. Under standard parabolic well-posedness assumptions, viscous regularization yields smooth solutions for positive times, making the strong form suitable for approximation with smooth neural networks. However, as $\nu$ decreases, the solution develops extremely narrow internal layers that approximate discontinuities, which in turn makes training standard PINNs increasingly challenging in practice. Spectral bias and the associated optimization difficulties are central to this small-viscosity training problem~\cite{xu2019frequency, wang2022ntk}.

A promising approach to address this training challenge lies in lifting-based strategies (LBS), which represent a discontinuous solution $u(\mathbf{x})$ as the restriction of a smoother function $U(\mathbf{x},\mathbf{z})$ to an embedded graph $\mathcal{M} = \{(\mathbf{x}, \mathbf{z}(\mathbf{x}))\}$. Here, the auxiliary variable $\mathbf{z}$ encodes interface-related information---such as shock locations or region identifiers---enabling a lifted reformulation of the original problem. Such representations have been used to capture solution jumps and derivative discontinuities in elliptic interface problems~\cite{hu2022discontinuity, tseng2023cusp}. This separates the sharp variation carried by the auxiliary variables from the smoother function represented by the network.

For scalar hyperbolic equations, the lift-and-embed framework~\cite{sun2024lift} constructs the embedded graph and evaluates residuals on it, with interface locations either specified or learned through trainable parameters. These works demonstrate how interface-aware auxiliary variables improve neural representations. For complex or dynamically evolving shock patterns, the remaining challenge is to construct such variables automatically without explicitly parameterizing the discontinuity interfaces.

Although originally designed for handling discontinuities, lifting-based methods can be naturally extended to problems involving steep yet continuous transition layers--commonly encountered in viscous regularizations of hyperbolic systems. To bridge the aforementioned gap and further advance the applicability of lifting strategies, we introduce a novel Two-stage Adaptive Lifting PINN (TAL-PINN).

The core idea of our approach is to augment the physical coordinates \((t,\mathbf{x})\) with an auxiliary variable \(\mathbf{z}(t,\mathbf{x})\) that is automatically constructed from geometry-aware \(r\)-adaptive coordinates. This design is inspired by traditional \(r\)-adaptive numerical methods~\cite{huang1994moving, huang2011adaptive, tang2003adaptive}, where coordinate transformations \(\boldsymbol{\xi}(t,\mathbf{x})\) redistribute mesh points according to local solution features, concentrating resolution near steep gradients or discontinuities. Such adaptive coordinates vary rapidly across sharp transition layers while remaining smooth elsewhere, effectively playing the same role as manually designed indicators in conventional lifting strategies.

We thus define \(\mathbf{z}:=\boldsymbol{\xi}(t,\mathbf{x})\), leading to a lifted representation of the solution as \(u(t,\mathbf{x}) = U\!(t,\mathbf{x}, \mathbf{z}(t,\mathbf{x}))\). The associated embedding graph \(\mathcal{M} = \{(t,\mathbf{x}, \mathbf{z}(t,\mathbf{x}))\}\) smoothly encodes sharp features in a higher-dimensional space. Existing adaptive sampling strategies use residual or gradient information to concentrate training points near difficult regions~\cite{mao2023pinns_adaptive, hanna2022residual, wu2023comprehensive}. Here, we parameterize the embedded graph by the computational coordinates and sample those coordinates uniformly, producing an adaptive distribution in the physical domain without explicit residual tracking. The same coordinate field therefore supplies the lifting input and the adaptive sampling distribution, coupling representation learning with geometry-induced point allocation.

A practical challenge arises because the lifting field depends on the unknown solution. TAL-PINN addresses this through two training stages. Stage 1 trains the solution network with identity coordinates at a relatively large viscosity, yielding a coarse solution whose broader transition layers are easier to resolve. This solution provides the geometric information used to construct the adaptive mesh and fit the coordinate network. Stage 2 uses the trained solution-network parameters as its initialization and the fitted coordinates as lifting inputs for training at the target small viscosity.

This progression combines the curriculum-learning principle of moving from easier to harder training problems~\cite{krishnapriyan2021failure} with solution-dependent lifting. It extends naturally to viscosity continuation, in which the viscosity decreases over successive training stages and the solution-network parameters are carried forward. In TAL-PINN, adaptive coordinates are also updated between stages.

We derive a scalar a posteriori error estimate relating the network error to the PDE residual and viscosity, and a variance-reduction criterion for adaptive residual sampling. A neural tangent kernel (NTK) analysis of the gradient-flow dynamics then explains how the lifting coordinates change the tangent features and the decay of residual modes.
Numerical experiments on the 1D Burgers equation with a stationary shock examine the sampling and optimization mechanisms through NTK spectra and loss decay, with ablations identifying the contribution of lifting beyond adaptive sampling and viscosity continuation. Additional experiments on the 1D moving-shock Burgers equation, the 1D Euler Sod and Lax shock tubes, and the 2D scalar Burgers equation demonstrate the practical effectiveness of TAL-PINN by accurately approximating near-inviscid solutions and reconstructing shock structures.

The remainder of this paper is structured as follows. Section~\ref{sec:preliminaries} introduces the problem background and key preliminaries. Section~\ref{sec:proposedmethod} develops the TAL-PINN framework, beginning with the general lifted formulation and residual, followed by the construction of $r$-adaptive coordinates, their use for lifting and sampling, and the two-stage training scheme with its multistage extension. Section~\ref{sec:theoretical} presents the theoretical analysis, providing an a posteriori error estimate, an importance-sampling interpretation, and an NTK-based gradient-flow analysis. Numerical experiments in Section~\ref{sec:numerical_experiments} validate the method's efficacy, and Section~\ref{sec:conclusion} offers concluding remarks.

\section{Preliminaries}
\label{sec:preliminaries}

In this section, we introduce viscous approximations of hyperbolic conservation laws, the PINN formulation, and lifting representations that form the basis of the proposed method.

\subsection{Hyperbolic conservation laws and viscous regularization}\label{sec:viscous_prelim}

We consider hyperbolic systems of conservation laws in the form
\begin{equation}
    \partial_t u + \nabla \cdot f(u) = 0,
    \label{eq:conservation}
\end{equation}
where \( u \in \mathbb{R}^m \) represents the vector of conserved variables, and \( f: \mathbb{R}^m \to \mathbb{R}^{m \times d} \) is the associated flux tensor~\cite{dafermos2005hyperbolic, smoller2012shock}.
Even with smooth initial data, classical solutions typically develop discontinuities in finite time. Beyond the formation of shocks, the notion of solutions is extended to weak solutions, which, however, are generally non-unique. An admissibility criterion---most commonly the entropy condition---rules out non-physical weak solutions. We denote the entropy-admissible solution of \eqref{eq:conservation} considered here by \(u^*\).

A widely adopted mechanism to enforce entropy admissibility is through viscous regularization. Artificial viscosity has also long been used for numerical shock capturing~\cite{vonneumann1950method}. Here we consider
\begin{equation}
    \partial_t u + \nabla \cdot f(u) = \nu \Delta u, \quad \nu > 0,
    \label{eq:viscous_conservation}
\end{equation}
where \( \nu \) is the viscosity coefficient. Under standard well-posedness assumptions, let \(u^\nu\) denote the viscous solution of \eqref{eq:viscous_conservation}, smooth for positive times. For scalar conservation laws and systems satisfying the corresponding vanishing-viscosity hypotheses, \(u^\nu \to u^*\) as \(\nu \to 0^+\) in the appropriate topology, selecting an entropy-admissible inviscid solution~\cite{bianchini2005vanishing}.

When the viscosity is small, the solution develops narrow internal layers with steep gradients that closely approximate discontinuities. Although these layers remain continuous, they pose significant numerical challenges. Traditional discretization methods, such as finite difference and finite element schemes, require extremely fine meshes to resolve such features. In high-dimensional settings, maintaining this resolution exacerbates the curse of dimensionality, leading to rapid growth in computational and memory demands~\cite{hirsch2007numerical, huang2011adaptive}.

\subsection{Physics-Informed Neural Networks (PINNs)}

In the PINN framework, the viscous solution \(u^\nu(t,\mathbf{x})\) of \eqref{eq:viscous_conservation} is approximated by a neural network \( \tilde{u}(t, \mathbf{x}; \theta) \), parameterized by \(\theta\). A fully connected feedforward neural network is typically employed:
\begin{equation}
    \tilde{u}(t, \mathbf{x}; \theta) = W_L \sigma \big( W_{L-1} \sigma ( \cdots \sigma ( W_1 [t, \mathbf{x}]^\top + b_1 ) \cdots ) + b_{L-1} \big) + b_L,
    \label{eq:pinn_architecture}
\end{equation}
where \( W_j \) and \( b_j \) denote the weights and biases of the \(j\)-th layer, \( L \) is the total number of layers, and \( \sigma(\cdot) \) is a nonlinear activation function.

The governing PDE is enforced through a residual formulation:
\begin{equation}
    \mathcal{R}(t, \mathbf{x}; \theta) := \partial_t \tilde{u}(t, \mathbf{x}; \theta) + \nabla \cdot f(\tilde{u}(t, \mathbf{x}; \theta)) - \nu\Delta_{\mathbf{x}}\tilde{u}(t,\mathbf{x};\theta),
    \label{eq:pde_residual}
\end{equation}
where all derivatives are computed using automatic differentiation.

We use $\ell$ for empirical losses evaluated on finite samples and $\mathcal{L}$ for their continuous counterparts. The empirical training loss consists of three components: the PDE residual, the initial condition (IC), and the boundary condition (BC):
\begin{equation}
    \ell(\theta)
    = \ell_{r}(\theta)
    + w_{\mathrm{ic}}\,\ell_{\mathrm{ic}}(\theta)
    + w_{\mathrm{bc}}\,\ell_{\mathrm{bc}}(\theta),
    \label{eq:pinn_loss_three}
\end{equation}
where \(w_{\mathrm{ic}}\) and \(w_{\mathrm{bc}}\) balance the IC and BC terms.
Each term is evaluated on its own set of collocation points:
\begin{equation}
\begin{aligned}
    \ell_{r}(\theta)
    &= \frac{1}{N_{r}} \sum_{i=1}^{N_{r}}
       \left| \mathcal{R}\!\left(t_i^{\,r}, \mathbf{x}_i^{\,r}; \theta\right) \right|^{2},\\[3pt]
    \ell_{\mathrm{ic}}(\theta)
    &= \frac{1}{N_{\mathrm{ic}}} \sum_{i=1}^{N_{\mathrm{ic}}}
       \left| \tilde{u}\!\left(0, \mathbf{x}_i^{\,\mathrm{ic}}; \theta\right) - u_{0}\!\left(\mathbf{x}_i^{\,\mathrm{ic}}\right) \right|^{2},\\[3pt]
    \ell_{\mathrm{bc}}(\theta)
    &= \frac{1}{N_{\mathrm{bc}}} \sum_{i=1}^{N_{\mathrm{bc}}}
       \left| \tilde{u}\!\left(t_i^{\,\mathrm{bc}}, \mathbf{x}_i^{\,\mathrm{bc}}; \theta\right) - g_{\mathrm{bc}}\!\left(t_i^{\,\mathrm{bc}}, \mathbf{x}_i^{\,\mathrm{bc}}\right) \right|^{2},
\end{aligned}
\label{eq:pinn_loss_terms_three}
\end{equation}
where \(\{(t_i^{\,r},\mathbf{x}_i^{\,r})\}_{i=1}^{N_r}\), \(\{(0,\mathbf{x}_i^{\,\mathrm{ic}})\}_{i=1}^{N_{\mathrm{ic}}}\), and \(\{(t_i^{\,\mathrm{bc}},\mathbf{x}_i^{\,\mathrm{bc}})\}_{i=1}^{N_{\mathrm{bc}}}\) are the residual, IC, and BC collocation sets, respectively; \(u_0(\cdot)\) is the prescribed initial data; and \(g_{\mathrm{bc}}(t,\mathbf{x})\) encodes the boundary data.

At small viscosity, the narrow transition layers introduce high-frequency components that are difficult for networks to learn. The Frequency Principle (F-Principle) describes the tendency of neural networks to learn low-frequency components faster than high-frequency ones~\cite{xu2019frequency, xu2024overview, xu2025understanding}. This motivates lifting representations in which auxiliary variables encode the sharp variation across these layers.

\subsection{Lifting-based strategies}

A common strategy for handling discontinuities is to represent the solution in an augmented space. Specifically, a possibly discontinuous function $u(\mathbf{x})$ is expressed by evaluating a smooth function $U$ on the lifted graph $\mathcal{M}=\{(\mathbf{x},\mathbf{z}(\mathbf{x}))\}$:
\begin{equation}
    u(\mathbf{x})=U(\mathbf{x},\mathbf{z}(\mathbf{x})),
    \label{eq:lift_formulation}
\end{equation}
where the auxiliary variable $\mathbf{z}(\mathbf{x})\in\mathbb{R}^k$ encodes local information about a shock or interface.

\begin{figure}[htbp]
    \centering
    \begin{subfigure}[b]{0.480\linewidth}
        \centering
        \includegraphics[width=\linewidth]{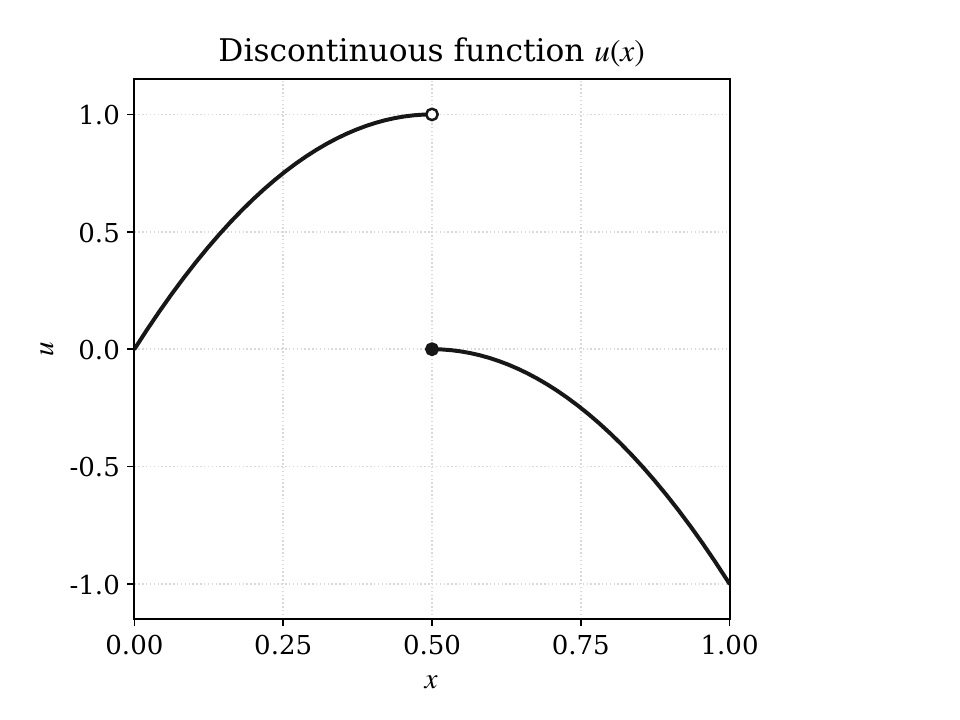}
        \caption{Discontinuous function $u(x)$.}
    \end{subfigure}
    \hfill
    \begin{subfigure}[b]{0.480\linewidth}
        \centering
        \includegraphics[width=\linewidth]{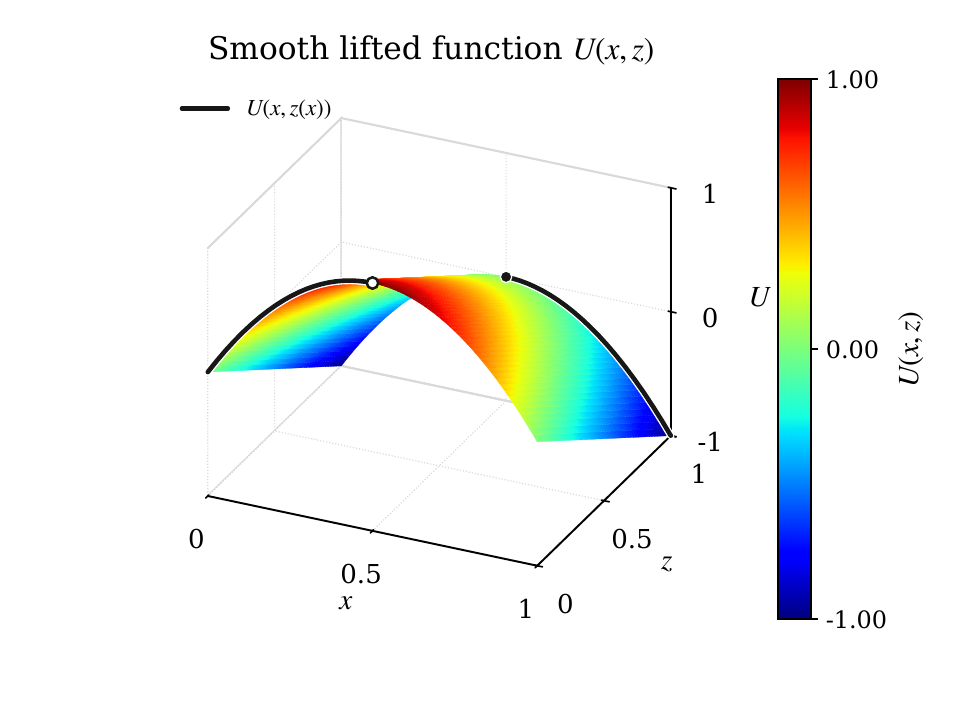}
        \caption{Smooth lifted function $U(x,z)$.}
    \end{subfigure}
    \caption{Lifting representation of $u(x)=4x(1-x)-H(x-0.5)$. The colored surface is $U(x,z)=4x(1-x)-z$, and the black curves show its restriction to $z=H(x-0.5)$. The open and filled circles indicate the left limit and the attained value at $x=0.5$, respectively. The two segments of the lifted graph are disconnected.}
    \label{fig:liftdemo}
\end{figure}

Figure~\ref{fig:liftdemo} illustrates this idea with the piecewise smooth function
\[
    u(x)=4x(1-x)-H(x-0.5),\qquad x\in[0,1],
\]
where $H(s)=0$ for $s<0$ and $H(s)=1$ for $s\geq0$. Introducing $z(x)=H(x-0.5)$ gives
\[
    u(x)=U(x,z(x)),\qquad U(x,z)=4x(1-x)-z.
\]
The function $U$ is smooth throughout the extended $(x,z)$-domain. The two branches of $u$ are recovered by evaluating $U$ at $z=0$ for $x<0.5$ and at $z=1$ for $x\geq0.5$. Thus the discontinuity of the composed function is encoded in the auxiliary variable, while the function represented in the augmented space remains smooth.

Where the auxiliary mapping is differentiable, the gradient of the composed function is obtained by the chain rule:
\begin{equation}
    \nabla_{\mathbf{x}}u
    =\nabla_{\mathbf{x}}U
    +(\nabla_{\mathbf{x}}\mathbf{z})^{\top}\nabla_{\mathbf{z}}U,
    \label{eq:lift_chainrule}
\end{equation}
where $\nabla_{\mathbf{x}}\mathbf{z}\in\mathbb{R}^{k\times d}$ is the Jacobian of the auxiliary mapping, with entries $(\nabla_{\mathbf{x}}\mathbf{z})_{\alpha j}=\partial_{x_j}z_\alpha$, and the derivatives of $U$ are evaluated at $(\mathbf{x},\mathbf{z}(\mathbf{x}))$. For the step-based example above, this identity applies separately on each side of the interface. In the viscously regularized setting considered below, smooth auxiliary coordinates allow the required derivatives to be evaluated throughout the training domain. Higher-order derivatives follow similarly, yielding a lifted residual while retaining the original differential operator. Related lifting representations have been investigated in~\cite{hu2022discontinuity,tseng2023cusp,sun2024lift}.

\section{Proposed method}
\label{sec:proposedmethod}

This section is structured to progressively develop the theoretical and algorithmic components of the proposed TAL-PINN framework. We first formulate the viscously regularized problem using a general smooth lifting variable. We then review $r$-adaptive coordinate construction and explain how these coordinates supply both the lifting input and the sampling distribution. Finally, we assemble these components into the two-stage training procedure and its multistage extension.

\subsection{Combining the lifting representation and viscous regularization}
\label{sec:lifted_formulation}

The lifting representation encodes local discontinuity information into an auxiliary variable $\mathbf{z}(t, \mathbf{x})$, yielding a continuous representation of the solution in the augmented space. We approximate this lifted solution using a solution network $\tilde{U}(t, \mathbf{x}, \mathbf{z}; \theta)$ with trainable parameters $\theta$. In viscous regularizations of hyperbolic conservation laws with small viscosity coefficients, singularities present in the inviscid limit manifest as narrow internal layers with steep gradients. Building on this structure, we propose to embed such high-gradient features directly into the auxiliary variable $\mathbf{z}(t, \mathbf{x})$, which promotes smoothness of the lifted solution $\tilde{U}(t, \mathbf{x}, \mathbf{z}; \theta)$ in the augmented space.

To illustrate how a smooth lifting variable encodes a viscous layer, consider a shock interface locally described by $x_d=\gamma(t,x_1,\dots,x_{d-1})$. The function $s(t,\mathbf{x}):=x_d-\gamma(t,x_1,\dots,x_{d-1})$ locates the transition region.
For a nondegenerate shock of fixed strength, adding a viscous regularization term \( \nu \Delta u \) produces a transition layer of width \( O(\nu) \) near the shock location. 
To capture this high-gradient structure, we introduce an auxiliary variable and define a lifted representation $\tilde{U}(t, \mathbf{x}, z(t, \mathbf{x}); \theta)$. For a scalar auxiliary coordinate, we take $z(t,\mathbf{x}) = \psi_\epsilon(s(t,\mathbf{x}))$, where $\psi_\epsilon: \mathbb{R} \to (0,1)$ is the smooth transition function
\begin{equation*}
\psi_\epsilon(s) = \frac{1}{2}\left[1+\tanh\!\left(\frac{s}{\epsilon}\right)\right],\qquad \epsilon>0.
\end{equation*}
This function is $C^\infty$ and encodes the local transition into the auxiliary coordinate. The parameter \(\epsilon\) controls the width of the encoded transition region and can be chosen proportional to the viscous-layer thickness, i.e., \(\epsilon = O(\nu)\). For a smooth interface function $s$, this construction supplies the derivatives required by the viscous residual below.

More generally, for a given smooth auxiliary field $\mathbf{z}(t,\mathbf{x})\in\mathbb{R}^k$, define the embedding
\[
    \Phi:(t,\mathbf{x})\mapsto(t,\mathbf{x},\mathbf{z}(t,\mathbf{x})).
\]
Evaluating the solution network on the lifted graph $\mathcal{M}=\{(t,\mathbf{x},\mathbf{z}(t,\mathbf{x}))\}$ gives the physical approximation
\begin{equation*}
\tilde{u}(t,\mathbf{x};\theta)
= \big(\tilde{U}(\cdot;\theta)\circ\Phi\big)(t,\mathbf{x})
= \tilde{U}(t,\mathbf{x},\mathbf{z}(t,\mathbf{x});\theta).
\end{equation*}

Using the convention that \( \nabla_{\mathbf{x}}\mathbf{z} \in \mathbb{R}^{k\times d} \) stacks the \(\mathbf{x}\)-gradients of each component of \(\mathbf{z}\) as rows, the chain rule gives
\begin{equation*}
\begin{aligned}
    \partial_t \tilde{u} &= \partial_t \tilde{U}
    + (\partial_t \mathbf{z})^\top \nabla_{\mathbf{z}} \tilde{U}, \quad
    \nabla_{\mathbf{x}} \tilde{u} = \nabla_{\mathbf{x}} \tilde{U}
    + (\nabla_{\mathbf{x}}\mathbf{z})^\top \nabla_{\mathbf{z}} \tilde{U}, \\
    \Delta_{\mathbf{x}} \tilde{u} &= \Delta_{\mathbf{x}} \tilde{U}
    + 2\,\mathrm{tr}\!\big((\nabla_{\mathbf{x}}\mathbf{z})\nabla_{\mathbf{z}}\nabla_{\mathbf{x}}\tilde{U}\big)
    + \mathrm{tr}\!\big((\nabla_{\mathbf{x}}\mathbf{z})^\top \nabla_{\mathbf{z}}^2 \tilde{U}\, (\nabla_{\mathbf{x}}\mathbf{z})\big)+ (\Delta_{\mathbf{x}}\mathbf{z})^\top \nabla_{\mathbf{z}} \tilde{U},
\end{aligned}
\end{equation*}
where \( \nabla_{\mathbf{z}}\nabla_{\mathbf{x}}\tilde{U} \in \mathbb{R}^{d\times k} \) has entries {\( \partial^2\tilde{U}/(\partial x_i\partial z_\alpha) \)}, and \( \nabla_{\mathbf{z}}^2\tilde{U} \in \mathbb{R}^{k\times k} \) is the Hessian of \( \tilde{U} \) with respect to \( \mathbf{z} \). All derivatives are evaluated at \( (t,\mathbf{x},\mathbf{z}(t,\mathbf{x})) \). For vector-valued outputs \( \tilde{U}\in\mathbb{R}^m \), the above relations are applied componentwise.

We then define the lifted PDE residual associated with the solution network as
\begin{equation}\label{eq:resfunctionofU}
\begin{aligned}
\mathcal{R}(t,\mathbf{x};\theta)
&:= \partial_t \tilde{U}
  + (\partial_t \mathbf{z})^\top \nabla_{\mathbf{z}} \tilde{U}\\
&\quad+\sum_{j=1}^{d} Df_j(\tilde{U})
  \left(\partial_{x_j}\tilde{U}
  +\sum_{\alpha=1}^{k}(\partial_{x_j}z_\alpha)\partial_{z_\alpha}\tilde{U}\right)\\
&\quad-\nu\Big[
  \Delta_{\mathbf{x}}\tilde{U}
  +2\mathrm{tr}\!\big((\nabla_{\mathbf{x}}\mathbf{z})
             \nabla_{\mathbf{z}}\nabla_{\mathbf{x}}\tilde{U}\big)\\
&\qquad\quad
  +\mathrm{tr}\!\big((\nabla_{\mathbf{x}}\mathbf{z})^\top
        \nabla_{\mathbf{z}}^2\tilde{U}(\nabla_{\mathbf{x}}\mathbf{z})\big)
  +(\Delta_{\mathbf{x}}\mathbf{z})^\top\nabla_{\mathbf{z}}\tilde{U}
  \Big].
\end{aligned}
\end{equation}
Here $Df_j$ denotes the Jacobian of the flux in the $x_j$ direction; for systems, the time and diffusion terms are evaluated componentwise.

Similar to  \eqref{eq:pinn_loss_three}, the total loss is defined as a weighted sum of the PDE residual and the initial-boundary condition terms: 
\begin{equation}
\ell(\theta) = \ell_{r}(\theta) + w_{\mathrm{ic}}\,\ell_{\mathrm{ic}}(\theta)  + w_{\mathrm{bc}}\,\ell_{\mathrm{bc}}(\theta),
\label{eq:total_loss_act}
\end{equation}
where
\begin{align}
\ell_{r}(\theta)
&= \frac{1}{N_r}\sum_{i=1}^{N_r}
   w\!\left(t_i^{\,r},\mathbf{x}_i^{\,r}\right)
   \left|\mathcal{R}\!\left(t_i^{\,r}, \mathbf{x}_i^{\,r};\theta\right)\right|^{2},
\label{eq:loss_termofU}\\
\ell_{\mathrm{ic}}(\theta)
&= \frac{1}{N_{\mathrm{ic}}} \sum_{i=1}^{N_{\mathrm{ic}}}
   \left| \tilde{U}\!\left(0, \mathbf{x}_i^{\,\mathrm{ic}}, \mathbf{z}(0, \mathbf{x}_i^{\,\mathrm{ic}})  ; \theta\right)
   - u_{0}\!\left(\mathbf{x}_i^{\,\mathrm{ic}}\right) \right|^{2},
\label{eq:loss_termofU2}\\
\ell_{\mathrm{bc}}(\theta)
&= \frac{1}{N_{\mathrm{bc}}} \sum_{i=1}^{N_{\mathrm{bc}}}
   \left| \tilde{U}\!\left(t_i^{\,\mathrm{bc}}, \mathbf{x}_i^{\,\mathrm{bc}},  \mathbf{z}(t_i^{\,\mathrm{bc}}, \mathbf{x}_i^{\,\mathrm{bc}}); \theta\right)
   - g_{\mathrm{bc}}\!\left(t_i^{\,\mathrm{bc}}, \mathbf{x}_i^{\,\mathrm{bc}}\right) \right|^{2},
\label{eq:loss_termofU3}
\end{align}
and $w(t,\mathbf{x})>0$ accounts for the sampling distribution of the residual points.

Manually prescribed lifting variables such as $H(s(t,\mathbf{x}))$ and $\psi_\epsilon(s(t,\mathbf{x}))$ use an interface function $s$ specified in advance. For unknown or evolving interfaces, we instead construct the lifting field from the transition geometry in an approximate solution.

\subsection{$r$-adaptivity and the equidistribution principle}
\label{sec:r_adaptivity}

We recall the essentials of \( r \)-adaptivity needed for our lifting construction. Let the physical domain be \(\Omega_P \subset \mathbb{R}^d\) and a fixed computational domain $\Omega_C \subset \mathbb{R}^d$ satisfy $|\Omega_C| = 1$, with coordinates $\mathbf{x} \in \Omega_P$ and $\boldsymbol{\xi} \in \Omega_C$, respectively. An \( r \)-adaptive method employs a time-dependent diffeomorphic mapping
\begin{equation*}
\mathbf{x} = \mathbf{x}(t, \boldsymbol{\xi}): [0,T] \times \Omega_C \to \Omega_P,
\end{equation*}
which transforms a fixed mesh in $\Omega_C$ into a nonuniform mesh in $\Omega_P$. Its (smooth) inverse
\begin{equation}
\boldsymbol{\xi}=\boldsymbol{\xi}(t,\mathbf{x}):[0,T]\times\Omega_P\to\Omega_C
\label{eq:coor_inversemap}
\end{equation}
maps the physical domain back to the computational one.

Let $M(t,\mathbf{x})>0$ be a monitor function that becomes large near steep features of the solution. A common gradient-based choice is
\begin{equation*}
M(t, \mathbf{x}) = \sqrt{1 + \beta \, | \nabla_{\mathbf{x}} \tilde{u}(t, \mathbf{x}; \theta) |^2},
\qquad \beta > 0.
\end{equation*}

The equidistribution principle requires that the computational and physical ``mass'' of $M$ match under the mapping. This condition reads as, for any measurable subset $A\subset\Omega_C$,
\begin{equation*}
\displaystyle\frac{\int_A \rd\boldsymbol{\xi}}{\int_{\Omega_C} \rd\boldsymbol{\xi}}
\;=\;
\frac{\int_{\mathbf{x}(t,A)} M(t, \mathbf{x}) \rd\mathbf{x}}{\int_{\Omega_P} M(t, \mathbf{x}) \rd\mathbf{x}}.
\end{equation*}
By a change of variables, this condition becomes the point-wise relation
\begin{equation}
M\big(t, \mathbf{x}(t, \boldsymbol{\xi})\big)\, \mathcal{J}(t, \boldsymbol{\xi}) = \Theta(t),
\quad
\Theta(t) = \int_{\Omega_P} M(t, \mathbf{x}) \rd\mathbf{x},
\label{eq:equid}
\end{equation}
where \( \mathcal{J}(t, \boldsymbol{\xi}) := \det\!\big(\nabla_{\boldsymbol{\xi}}\mathbf{x}\big)>0 \) is the Jacobian determinant of the orientation-preserving transformation.

In one spatial dimension, let \( \Omega_P = [x_a, x_b] \) and \( \Omega_C = [0, 1] \). Strict equidistribution gives an explicit map:
\begin{equation}
\xi(t, x) = \frac{\int_{x_a}^{x} M(t, s) \ds}{\int_{x_a}^{x_b} M(t, s) \ds},
\qquad
x(t, \xi) = \big(\xi(t,\cdot)\big)^{-1}(\xi).
\label{eq:equidistribution_1d_map}
\end{equation}

In multiple dimensions, equidistribution alone does not determine a unique mapping. A regularizing elliptic surrogate balances concentration according to the monitor with mesh smoothness. A common choice is the isotropic Winslow-type variable-diffusion equation, posed at each fixed time $t$:
\begin{equation}
\nabla_{\boldsymbol{\xi}} \cdot \left( M\big(t, \mathbf{x}(t, \boldsymbol{\xi})\big)\, \nabla_{\boldsymbol{\xi}} \mathbf{x}(t, \boldsymbol{\xi}) \right) = 0
\quad \text{in } \Omega_C,
\quad
\mathbf{x}(t, \boldsymbol{\xi}) = \mathbf{x}_b(\boldsymbol{\xi}) \ \text{on } \partial \Omega_C.
\label{eq:winslow}
\end{equation}

A corresponding moving mesh partial differential equation (MMPDE) is obtained by relaxing \eqref{eq:winslow} in time. This MMPDE5-type evolution uses the elliptic condition to guide mesh redistribution and produces a discrete adaptive mesh~\cite{huang1994mmpde,huang2011adaptive}. The uniformly distributed computational nodes $\{\boldsymbol{\xi}_i\} \subset \Omega_C$ are mapped to physical nodes $\mathbf{x}_i = \mathbf{x}(t, \boldsymbol{\xi}_i)$ in $\Omega_P$. This yields a discrete representation $\{(\boldsymbol{\xi}_i, \mathbf{x}_i)\}$ of the diffeomorphism, producing a nonuniform physical mesh aligned with the monitor function $M(t,\mathbf{x})$.

\subsection{Adaptive lifting and coordinate-induced sampling}
\label{sec:lifting_variable}

The inverse map $\boldsymbol{\xi}(t,\mathbf{x})$ in \eqref{eq:coor_inversemap}, obtained from the adaptive mesh construction, varies gradually in regular regions and rapidly across shock-aligned transition layers. It plays a role similar to manually constructed smooth indicators such as $\psi_\epsilon$, with its transition determined by the solution-dependent monitor.

\begin{figure}[htbp]
    \centering
    \begin{subfigure}[b]{0.480\linewidth}
        \centering
        \includegraphics[width=\linewidth]{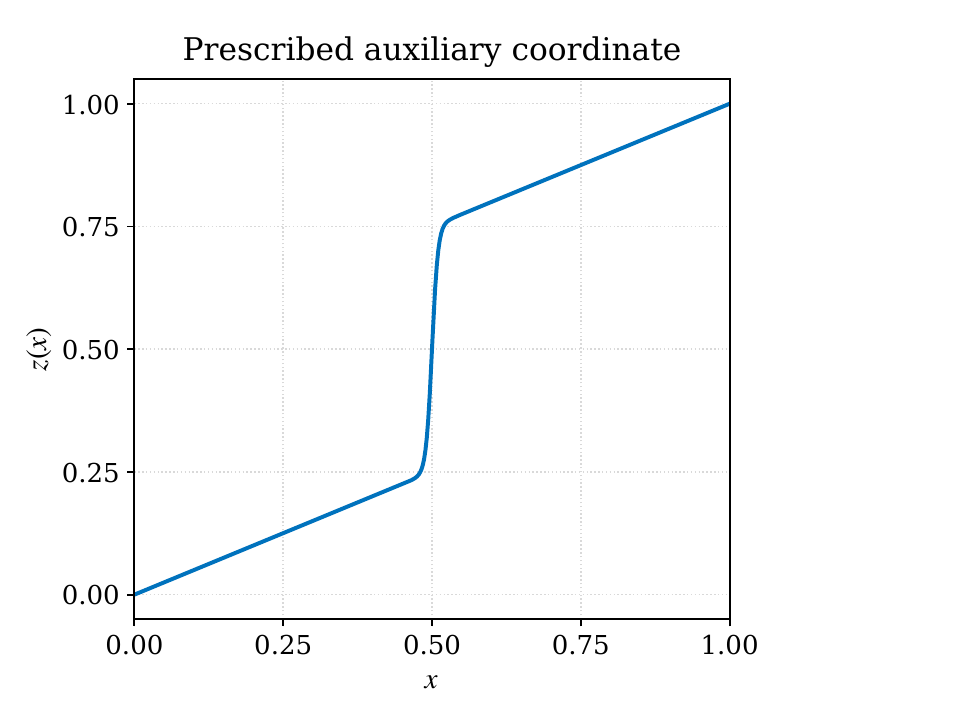}
    \end{subfigure}
    \hfill
    \begin{subfigure}[b]{0.480\linewidth}
        \centering
        \includegraphics[width=\linewidth]{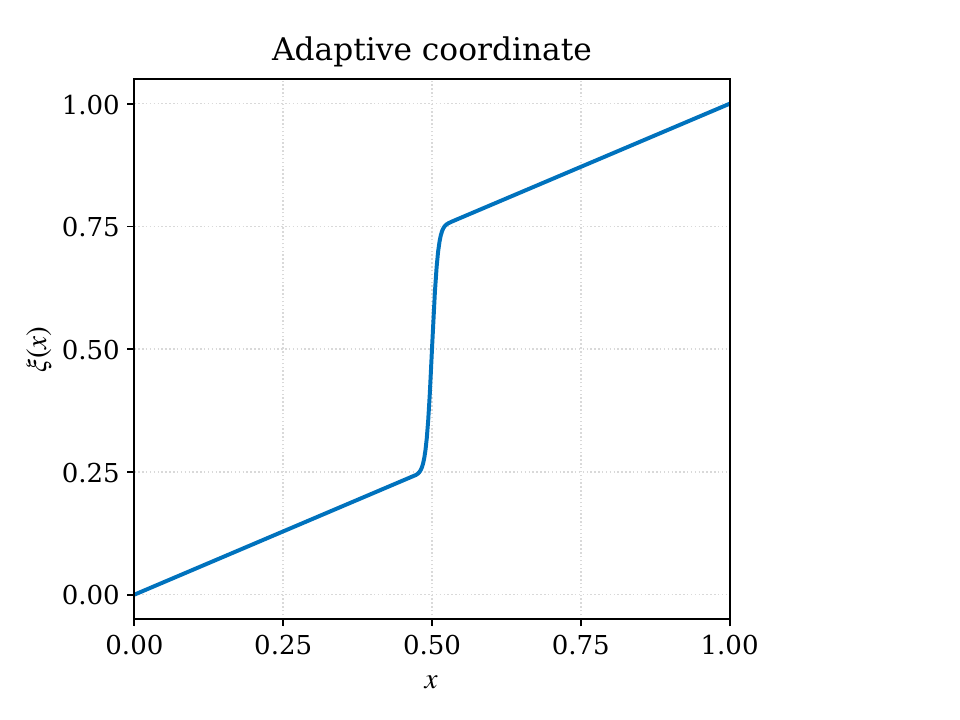}
    \end{subfigure}
    \caption{Left: prescribed auxiliary coordinate $z(x)=\tfrac12[x+\psi_\epsilon(x-x_0)]$, where $\psi_\epsilon(s)=\tfrac12[1+\tanh(s/\epsilon)]$, $x_0=0.5$ and $\epsilon=10^{-2}$. Right: adaptive coordinate $\xi(x)=\int_0^x M(s)\,ds/\int_0^1 M(s)\,ds$, with $M=\sqrt{1+u_x^2}$, constructed from the analytic stationary Burgers shock $u(x)=-a\tanh(a(x-x_0)/(2\nu))$, $a=0.5$ and $\nu=2.5\times10^{-3}$.}
    \label{fig:liftdemo2}
\end{figure}

For illustration in one spatial dimension, Figure~\ref{fig:liftdemo2} compares a prescribed auxiliary variable $z(x)=\tfrac{1}{2}\big(x+\psi_\epsilon(x-0.5)\big)$ with an adaptive coordinate $\xi(x)$ constructed from an analytic stationary Burgers shock. Both coordinates exhibit a narrow transition near the shock and vary gradually away from it. This similarity illustrates how the adaptive coordinate encodes the transition geometry through the monitor.

Motivated by this observation, we propose to use the adaptive coordinate itself as the lifting input: $
    \mathbf{z}(t, \mathbf{x}) := \boldsymbol{\xi}(t, \mathbf{x}) \in \mathbb{R}^d$,
which yields a lifted representation:
\begin{equation*}
    \tilde{u}(t, \mathbf{x}; \theta) 
    \;=\; \tilde{U}\big(t, \mathbf{x}, \boldsymbol{\xi}(t, \mathbf{x}); \theta\big),
\end{equation*}
implicitly encoding singular features of the solution.

Crucially, constructing \( \boldsymbol{\xi}(t, \mathbf{x}) \) does not require any a priori knowledge of the target small-viscosity solution. At a moderately larger viscosity $\nu_0$, the transition layers are wider and easier for the network to resolve, while retaining useful information about the locations and geometry of steep features. The initial high-viscosity training stage thus provides a coarse solution from which we construct the monitor function \( M(t, \mathbf{x}) \) and the associated adaptive coordinate \( \boldsymbol{\xi}(t, \mathbf{x}) \). This coordinate supplies a geometric prior through the lifting input for subsequent training in the small-viscosity regime.

The same adaptive coordinates also determine the sampling distribution of the residual points. With $\mathbf{z}=\boldsymbol{\xi}(t,\mathbf{x})$, the lifted graph at each fixed time admits the parameterization
\[
    \mathcal{M}_t
    =\{(\mathbf{x}(t,\boldsymbol{\xi}),\boldsymbol{\xi}):
          \boldsymbol{\xi}\in\Omega_C\}.
\]
For the smooth inverse maps introduced in Section~\ref{sec:r_adaptivity}, uniformly sampling $\boldsymbol{\xi}\in\Omega_C$ induces the physical density
\begin{equation}
    \rho_t(\mathbf{x})=\left|\det \nabla_{\mathbf{x}}
               \boldsymbol{\xi}(t,\mathbf{x})\right|,
    \qquad \int_{\Omega_P}\rho_t(\mathbf{x})\,\rd\mathbf{x}=1,
    \label{eq:coordinate_density}
\end{equation}
where we have used $|\Omega_C|=1$. The corresponding measure-correction weight is
\begin{equation}
    w(t,\mathbf{x})=\rho_t(\mathbf{x})^{-1}
    =\left|\det \nabla_{\mathbf{x}}\boldsymbol{\xi}(t,\mathbf{x})\right|^{-1}.
    \label{eq:measure_weight}
\end{equation}
This change of variables gives
\[
    \mathbb{E}_{\mathbf{x}\sim\rho_t}
       \big[w(t,\mathbf{x})|\mathcal{R}(t,\mathbf{x};\theta)|^2\big]
    =\int_{\Omega_P}|\mathcal{R}(t,\mathbf{x};\theta)|^2\,\rd\mathbf{x}.
\]
With time sampled uniformly on $[0,T]$, \eqref{eq:loss_termofU} approximates the time-averaged physical residual integral.

In practice, the computed \(r\)-adaptive mapping is available only at discrete nodes \( \{(t_i,\mathbf{x}_i,\boldsymbol{\xi}_i)\}_{i=1}^{N_{\xi}} \), where $N_{\xi}$ denotes the total number of discrete space-time nodes obtained from the moving mesh computation. Let $\mathbf{x}_h(t,\boldsymbol{\xi})$ denote the interpolant of the discrete forward-map data, providing a numerical approximation to $\mathbf{x}(t,\boldsymbol{\xi})$; the subscript $h$ indicates the mesh discretization. This map is used to transform uniform computational samples into physical residual points.

To obtain a smooth approximation of the inverse map $\boldsymbol{\xi}(t,\mathbf{x})$, we fit a lightweight coordinate neural network
\begin{equation}
\tilde{\boldsymbol{\xi}}(t,\mathbf{x};\theta_\xi):\mathbb{R}^{1+d}\to\mathbb{R}^{d}
\label{eq:coordnet}
\end{equation}
to the same nodal data, with $(t_i,\mathbf{x}_i)$ as inputs and $\boldsymbol{\xi}_i$ as targets. The fitted coordinate network supplies the lifting input and its derivatives at arbitrary query points for evaluating the lifted PDE residual. The associated training loss is defined as
\begin{equation}
\ell_{\xi}(\theta_\xi)
= \frac{1}{dN_{\xi}}\sum_{i=1}^{N_{\xi}}
\left\|\tilde{\boldsymbol{\xi}}(t_i,\mathbf{x}_i;\theta_\xi)-\boldsymbol{\xi}_i\right\|^2
+ \frac{w_J}{N_{\xi}}\sum_{i=1}^{N_{\xi}}
   \left[\max(0,\delta_J-\det J_{\xi,i})\right]^2,
\label{eq:loss_coordnet}
\end{equation}
where \(J_{\xi,i}:=\nabla_{\mathbf{x}}\tilde{\boldsymbol{\xi}}(t_i,\mathbf{x}_i;\theta_\xi)\in\mathbb{R}^{d\times d}\) is the spatial Jacobian and $\delta_J>0$. The coefficient \(w_J>0\) balances the mean-square fitting error against the anti-folding regularizer, which discourages small or negative Jacobian determinants; increasing $w_J$ more strongly penalizes mesh tangling.

Let \(\theta_{\xi}^*\) be an approximate minimizer of the loss function \(\ell_{\xi}(\theta_\xi)\). We define the lifting input as
$\mathbf{z}(t,\mathbf{x})=\tilde{\boldsymbol{\xi}}(t,\mathbf{x};\theta_{\xi}^*)$.
The required derivatives of $\mathbf{z}$ are then computed via  automatic differentiation:
\begin{equation}
\begin{aligned}
\partial_t \mathbf{z}(t,\mathbf{x})
&= \partial_t \tilde{\boldsymbol{\xi}}(t,\mathbf{x};\theta_{\xi}^*),\quad
\nabla_{\mathbf{x}} \mathbf{z}(t,\mathbf{x})
= \nabla_{\mathbf{x}} \tilde{\boldsymbol{\xi}}(t,\mathbf{x};\theta_{\xi}^*),\\
\Delta_{\mathbf{x}} \mathbf{z}(t,\mathbf{x})
&= \sum_{k=1}^{d}\partial_{x_kx_k}\tilde{\boldsymbol{\xi}}(t,\mathbf{x};\theta_{\xi}^*).
\end{aligned}
\label{eq:coordnet_derivs}
\end{equation}

The spatial Jacobian computed for the lifted residual also provides an approximation to the measure-correction weight in \eqref{eq:measure_weight}:
\begin{equation}
    \widehat w(t,\mathbf{x})
    =\left|\det \nabla_{\mathbf{x}}
       \tilde{\boldsymbol{\xi}}(t,\mathbf{x};\theta_\xi^*)\right|^{-1}.
    \label{eq:numerical_measure_weight}
\end{equation}
We apply these weights when assembling the residual loss in \eqref{eq:loss_termofU} from physical points generated by $\mathbf{x}_h$. For a domain-averaged loss, the corresponding self-normalized form is $\sum_i\widehat w_i|\mathcal{R}_i|^2/\sum_i\widehat w_i$. The fitted coordinates encode the transition geometry in the lifting input, while the adaptive mesh concentrates residual evaluations near steep layers. The weights account for the resulting sampling density in the physical residual objective.

\subsection{Two-stage training and multistage extension}
\label{sec:training_scheme}

We organize the adaptive lifting method into two training stages. In Stage 1, the solution network is trained at a coarse viscosity $\nu_0\gg\nu_{\mathrm{target}}$ with $\mathbf{z}(t,\mathbf{x})\equiv\mathbf{x}$ on the unit box, or its affine normalization on a general box. The resulting coarse solution defines the monitor and adaptive mesh used to fit the coordinate network through \eqref{eq:loss_coordnet}.

Stage 2 starts from the trained solution-network parameters and reduces the viscosity to $\nu_{\mathrm{target}}$. Its lifting input is $\mathbf{z}(t,\mathbf{x})=\tilde{\boldsymbol{\xi}}(t,\mathbf{x};\theta_\xi^*)$. The coordinate-network parameters $\theta_\xi^*$ remain fixed while the solution-network parameters $\theta$ are optimized using the lifted residual and loss in \eqref{eq:resfunctionofU}--\eqref{eq:loss_termofU3}.

Algorithm~\ref{alg:Lifting_framework} gives the two-stage procedure using a domain-averaged residual loss. Repeating the coordinate construction and lifted training yields a multistage viscosity continuation with $\nu_0>\nu_1>\cdots>\nu_K=\nu_{\mathrm{target}}$. Before training at each $\nu_k$, the preceding physical approximation is evaluated with its current lifting field to rebuild the monitor, adaptive mesh, and coordinate network. The solution-network parameters are carried over, while the optimizer state is restarted; the new coordinate-network parameters remain fixed throughout the next PDE-training stage.

\begin{algorithm}[H]
\small
\caption{TAL-PINN: two-stage adaptive lifting}
\label{alg:Lifting_framework}
\DontPrintSemicolon
\KwIn{PDE with initial/boundary conditions; viscosities $\nu_0\gg\nu_{\mathrm{target}}$; monitor prescription $M$; regular auxiliary grid $\{(t_j,\mathbf x_r^{a})\}$, uniform computational nodes $\{\boldsymbol\xi_i\}$, and training budgets.}
\KwOut{Solution network $\tilde U$ and lifting field $\mathbf z$.}
\textbf{Stage 1: coarse viscosity and coordinate construction}\;
Initialize $\theta$; set $\nu\leftarrow\nu_0$ and $\mathbf z_0(t,\mathbf x)=\mathbf x$ (with affine normalization)\;
\For{each Stage 1 training iteration}{
    Draw residual points uniformly on $[0,T]\times\Omega_P$ and initial/boundary points on their respective domains\;
    Evaluate \eqref{eq:resfunctionofU} with $\mathbf z=\mathbf z_0$; update $\theta$ using \eqref{eq:total_loss_act}--\eqref{eq:loss_termofU3} with residual weights $w_i=1$\;
}
Set $\theta_0^*\leftarrow\theta$ and evaluate the coarse solution on the auxiliary grid: $u_{jr}\leftarrow\tilde U(t_j,\mathbf x_r^{a},\mathbf z_0(t_j,\mathbf x_r^{a});\theta_0^*)$\;
Compute $M$ from the coarse solution; obtain physical nodes $\{\mathbf x_{ji}\}$ at $\{\boldsymbol\xi_i\}$ using \eqref{eq:equidistribution_1d_map} at each $t_j$ in 1D, or the MMPDE5-type evolution associated with \eqref{eq:winslow} in multiple dimensions\;
Interpolate $(t_j,\boldsymbol\xi_i)\mapsto\mathbf x_{ji}$ to obtain $\mathbf x_h(t,\boldsymbol\xi)$\;
Fit $\tilde{\boldsymbol\xi}(t,\mathbf x;\theta_\xi)$ to $\{(t_j,\mathbf x_{ji},\boldsymbol\xi_i)\}$ by minimizing \eqref{eq:loss_coordnet}; denote the fitted parameters by $\theta_\xi^*$\;
\BlankLine
\textbf{Stage 2: target viscosity with adaptive lifting}\;
Set $\theta\leftarrow\theta_0^*$, $\nu\leftarrow\nu_{\mathrm{target}}$, and $\mathbf z(t,\mathbf x)=\tilde{\boldsymbol\xi}(t,\mathbf x;\theta_\xi^*)$; restart the optimizer and keep $\theta_\xi^*$ fixed\;
\For{each Stage 2 training iteration}{
    Draw $t_i\sim\mathcal U([0,T])$ and $\boldsymbol\xi_i\sim\mathcal U(\Omega_C)$; set $\mathbf x_i\leftarrow\mathbf x_h(t_i,\boldsymbol\xi_i)$, $i=1,\ldots,N_r$\;
    Evaluate $\mathbf z_i=\mathbf z(t_i,\mathbf x_i)$ and its derivatives $\partial_t\mathbf z$, $\nabla_{\mathbf x}\mathbf z$, $\Delta_{\mathbf x}\mathbf z$ using \eqref{eq:coordnet_derivs}\;
    Compute $\mathcal R_i$ from \eqref{eq:resfunctionofU} at $(t_i,\mathbf x_i,\mathbf z_i)$ and $\widehat w_i=|\det\nabla_{\mathbf x}\mathbf z(t_i,\mathbf x_i)|^{-1}$\;
    Draw initial/boundary points as in Stage 1; evaluate $\ell_{\mathrm{ic}}$ and $\ell_{\mathrm{bc}}$ with the current lifting field\;
    Update only $\theta$ using $\nabla_\theta\ell$, where
    $\ell=\dfrac{\sum_i\widehat w_i|\mathcal R_i|^2}{\sum_i\widehat w_i}
    +w_{\mathrm{ic}}\ell_{\mathrm{ic}}+w_{\mathrm{bc}}\ell_{\mathrm{bc}}$\;
}
\KwRet{$\tilde U(\cdot;\theta)$ and $\mathbf z$, giving $\tilde u(t,\mathbf x)=\tilde U(t,\mathbf x,\mathbf z(t,\mathbf x);\theta)$}\;
\end{algorithm}

\section{Theoretical analysis}
\label{sec:theoretical}

In this section, we establish the theoretical foundations for the TAL-PINN framework presented in Section~\ref{sec:proposedmethod}. We first derive an a posteriori error estimate that rigorously quantifies how the network approximation error depends on both the PDE residual and the viscosity parameter. We further derive a variance-reduction criterion for adaptive sampling and explain how lifting changes the tangent-feature representation to accelerate residual decay. For simplicity, we write \(A \lesssim B\) to indicate that \(A \le C B\) for some constant \(C>0\), independent of the variables under consideration, unless otherwise specified.

\subsection{A posteriori error estimate}

For the scalar error analysis, let the physical domain \(\Omega_P\subset\mathbb{R}^d\) be bounded and Lipschitz, and let \(T>0\). Consider the viscously regularized equation on $\Omega_P\times(0,T)$:
\begin{equation}
    u^{\nu}_t + \nabla \cdot f(u^{\nu}) = \nu \Delta u^{\nu}
\label{eq:convection_diffusion_pde}
\end{equation}
with viscosity coefficient \(\nu > 0\) and a continuously differentiable flux \(f: \mathbb{R} \to \mathbb{R}^d\). We impose either homogeneous Dirichlet boundary conditions \((u^{\nu} = 0 \text{ on } \partial \Omega_P)\) or periodic boundary conditions.

Let \(u\) be an arbitrary approximation to the exact solution \(u^{\nu}\) of \eqref{eq:convection_diffusion_pde} that exactly satisfies the boundary conditions. In the periodic case, assume that $u$ and $u^\nu$ have the same spatial mean for every $t\in[0,T]$, so that $e=u^\nu-u$ has zero mean. Its residual is defined as
\[
    \mathcal{R}(u) := u_t + \nabla \cdot f(u) - \nu \Delta u.
\]

\begin{lemma}\label{thm:apost_L2}
Suppose $u \in L^2(0,T; H^1(\Omega_P))$,
$u_t \in L^2(0,T; H^{-1}(\Omega_P))$,
and $\mathcal{R}(u) \in L^2\big(0,T; L^2(\Omega_P)\big).
$
Assume also that $u$ and $u^\nu$ take values in a common bounded interval $I$. Then, for all $t\in[0,T]$, the error $e:=u^{\nu}-u$ satisfies
\[
\|e(t)\|_{L^2(\Omega_P)} \le
\exp\!\left(\frac{L_f^2}{2\nu} t\right)
\left[
\|e(0)\|_{L^2(\Omega_P)}
+\frac{C_P}{\sqrt{\nu}}\|\mathcal{R}(u)\|_{L^2(0,t;\,L^2(\Omega_P))}
\right],
\]
where \(C_P > 0\) is the Poincaré constant depending only on \(\Omega_P\) and the boundary condition, and $L_f := \sup_{s \in I} \left\| f'(s) \right\|_{\mathbb{R}^d}$. Here $\|\cdot\|_{\mathbb R^d}$ denotes the Euclidean norm.
\end{lemma}

\begin{proof}
From the equation of $u^{\nu}$ and the definition of $ \mathcal R(u)$, it follows that
\[
e_t+\nabla\!\cdot\!\big(f(u^{\nu})-f(u)\big)=\nu\Delta e - \mathcal R(u).
\]
Using the matching Dirichlet or periodic boundary conditions, testing the above equation with $e$ and applying integration by parts leads to
\[
\begin{split}
\frac12\frac{\rd}{\rd t}\|e\|_{L^2}^2+\nu\|\nabla e\|_{L^2}^2
&= \int_{\Omega_P}(f(u^{\nu})-f(u))\cdot\nabla e\rd\mathbf{x}
- \int_{\Omega_P}\mathcal R(u)e\rd\mathbf{x}\\
&\le L_f\|e\|_{L^2}\|\nabla e\|_{L^2}+C_P \|\mathcal R(u)\|_{L^2} \|\nabla e\|_{L^2}.
\end{split}
\]
Here we have used the Poincaré inequality. Further, applying Young's inequality with parameter $\nu$, we have
\[
\begin{aligned}
L_f\|e\|_{L^2}\|\nabla e\|_{L^2}
&\le \frac{\nu}{2}\|\nabla e\|_{L^2}^2+\frac{L_f^2}{2\nu}\|e\|_{L^2}^2,\\
C_P \|\mathcal R(u)\|_{L^2} \|\nabla e\|_{L^2}
&\le \frac{\nu}{2}\|\nabla e\|_{L^2}^2+\frac{C_P^2}{2\nu}\|\mathcal R(u)\|_{L^2}^2.
\end{aligned}
\]
Substituting and canceling $\nu\|\nabla e\|_{L^2}^2$ gives
\[
\frac12\frac{\rd}{\rd t}\|e\|_{L^2}^2
\le \frac{L_f^2}{2\nu}\|e\|_{L^2}^2+\frac{C_P^2}{2\nu}\|\mathcal R(u)\|_{L^2}^2.
\]
Applying Grönwall's inequality over \([0,t]\) yields
\[
\|e(t)\|_{L^2(\Omega_P)}^2 \le
\exp\!\left(\frac{L_f^2}{\nu} t\right)
\left[
\|e(0)\|_{L^2(\Omega_P)}^2
+\frac{C_P^2}{\nu}\int_0^t\|\mathcal R(u)(s)\|_{L^2(\Omega_P)}^2\,\ds
\right],
\]
which completes the proof.
\end{proof}

We conclude from Lemma~\ref{thm:apost_L2} that, for each fixed $\nu>0$, the PINN approximation converges to the viscous solution in $L^2(\Omega_P)$ uniformly over $[0,T]$ as the continuous initial and residual losses tend to zero. Combined with the classical theory of vanishing viscosity, this gives a route to approximating the inviscid entropy solution by controlling both the viscous approximation error and the network error.

This estimate also clarifies how viscosity governs error control. The bound contains both $1/\sqrt{\nu}$ and $\exp(L_f^2T/(2\nu))$, so its coefficients increase as $\nu\to0$. The estimate therefore requires stronger residual control to ensure the same accuracy at a smaller viscosity. This sensitivity motivates the continuation and lifting strategy for training near the inviscid limit.

Applying Lemma~\ref{thm:apost_L2} to the neural network approximation $\tilde{u}(t,\mathbf{x};\theta)$ in this setting and integrating over time gives a bound in terms of the continuous residual and initial losses:
\begin{equation*}
\|u^{\nu} - \tilde{u}(\cdot;\theta)\|^2_{L^2(0,T;\, L^2(\Omega_P))}
\le C_\nu\mathcal{L},\qquad \mathcal{L}:=\mathcal{L}_r+\mathcal{L}_{\mathrm{ic}},
\end{equation*}
where $C_\nu=T\exp(L_f^2T/\nu)\max(1,C_P^2/\nu)$ and
\begin{equation}
\begin{aligned}
\mathcal{L}_{r} &:= \int_0^T \int_{\Omega_P} \left| \mathcal{R}\!\big(\tilde{u}(t,\mathbf{x};\theta)\big) \right|^2  \rd\mathbf{x}\rd t, \\
\mathcal{L}_{\mathrm{ic}} &:= \int_{\Omega_P} \left| u^{\nu}(0,\mathbf{x}) - \tilde{u}(0,\mathbf{x};\theta) \right|^2  \rd\mathbf{x}.
\end{aligned}
\label{eq:total_continuous_loss}
\end{equation}
Further, its empirical counterpart $\ell_n:=\ell_{n,r}+\ell_{n,\mathrm{ic}}$ is defined by finite collocation samples:
\[
\begin{aligned}
\ell_{n,r}
&:=\frac{1}{N_r}\sum_{i=1}^{N_r}
\frac{\left|\mathcal{R}\big(\tilde u(t_i,\mathbf{x}_i;\theta)\big)\right|^2}
     {q_r(t_i,\mathbf{x}_i)},\\
\ell_{n,\mathrm{ic}}
&:=\frac{1}{N_{\mathrm{ic}}}\sum_{j=1}^{N_{\mathrm{ic}}}
\frac{\left|u^\nu(0,\mathbf{x}_j)-\tilde u(0,\mathbf{x}_j;\theta)\right|^2}
     {q_{\mathrm{ic}}(\mathbf{x}_j)}.
\end{aligned}
\]
Here $q_r$ and $q_{\mathrm{ic}}$ are the normalized sampling densities for the residual and initial points. Uniform time sampling with the adaptive spatial map gives $q_r(t,\mathbf{x})=\rho_t(\mathbf{x})/T$, while uniform initial sampling gives $q_{\mathrm{ic}}=1/|\Omega_P|$. These densities account for the time and domain factors that convert the averaged collocation losses into estimates of the integrals above.

Denote by \(u^*\) the exact entropy solution (inviscid case), \(u^\nu\) the viscous solution (\(\nu > 0\)), and \(\tilde{u}(t,\mathbf{x};\theta)\) the PINN approximation. The overall approximation error can be decomposed as
\begin{equation*}
\begin{aligned}
\|u^*-\tilde u(\cdot;\theta)\|^2_{L^2(0,T;\,L^2(\Omega_P))}
&\le 2\|u^*-u^\nu\|^2_{L^2(0,T;\,L^2(\Omega_P))}\\
&\quad+2C_\nu\left(\left|\mathcal{L}-\ell_n\right|+\ell_n\right).
\end{aligned}
\end{equation*}
Here, the first term is the vanishing-viscosity error, the second contribution is the continuous--empirical loss gap, and the third is the empirical training loss, with the latter two amplified by $C_\nu$. For uniformly bounded scalar solutions, $L^1$ vanishing-viscosity convergence also gives $L^2$ convergence through $\|u^\nu-u^*\|_{L^2}^2\le 2B_u\|u^\nu-u^*\|_{L^1}$, where $|u^\nu|,|u^*|\le B_u$. Thus, choosing the training accuracy so that $C_\nu(|\mathcal{L}-\ell_n|+\ell_n)\to0$ completes the error-control argument as $\nu\to0$.

\subsection{Importance sampling to reduce statistical error}
\label{sec:sec42}
In this subsection, we turn our attention to the PDE residual component of the loss function. Define the squared residual as
\[
g(y):=|\mathcal{R}(y)|^2,\quad\text{for}\, y=(t,\mathbf{x})\in D:=[0,T]\times\Omega_P,
\]
and let all sampling densities \(\rho\) be normalized, satisfying \(\displaystyle\int_D \rho\dy=1\). For simplicity, we assume \(|D|=1\). Consequently, the uniform density is \(\rho_u(y)\equiv1\).

We approximate the continuous loss
$
\displaystyle\mathcal{L}_r:=\int_{D} g(y)\dy
$
by the empirical estimator
\begin{equation*}
\ell_{n,r}:=\frac{1}{n}\sum_{i=1}^n \frac{g(y_i)}{\rho(y_i)},\quad y_i\overset{\mathrm{iid}}{\sim}\rho.
\end{equation*}
For fixed $g$ and $\rho$, assuming that \(\dfrac{g(y)}{\rho(y)}\) is bounded above by some constant \(B>0\), the Bernstein inequality yields~\cite{boucheron2003concentration} 
\begin{equation*}
\mathbb{P}\!\left(|\mathcal{L}_r-\ell_{n,r}|\ge \epsilon\right)
\ \le\ 2\exp\!\left(-\frac{n\epsilon^2}{2\sigma^2+\frac{2}{3}B\epsilon}\right),
\quad
\sigma^2:=\mathrm{Var}_{\rho}\!\left(\frac{g(y)}{\rho(y)}\right),\ \ y\sim\rho.
\end{equation*}
Consequently,  for a fixed sample size \(n\) and a common bound $B$, a smaller variance \(\sigma^2\) leads to a tighter probabilistic control over the estimation \(|\mathcal{L}_r-\ell_{n,r}|\).

In the vanilla PINN setting, collocation points are typically sampled from a uniform distribution, i.e. \(\rho(y)\equiv \rho_u(y)\). In this case, the variance becomes
\[
\mathrm{Var}_{\rho_u}(g(y))=\int_D g^2(y)\dy-\left(\int_D g(y)\dy\right)^2,
\]
which can attain large values when the residual \(g\) exhibits sharp local peaks or multiscale variabilities.

Adaptive sampling mitigates this issue through the principle of importance sampling, which aims to choose a density $\rho$ that correlates with $g$~\cite{robert2004monte}. Under the assumption that $\rho(y) > 0$ for almost every $y \in D$, the variance of the  importance-sampled estimator is
\begin{equation*}
\mathrm{Var}_{\rho}\!\left(\frac{g(y)}{\rho(y)}\right) = \int_D \frac{g^2(y)}{\rho(y)}\dy - \left(\int_D g(y)\dy\right)^2.
\end{equation*}

By the Cauchy--Schwarz inequality, we have
\begin{equation*}
\int_D \frac{g^2(y)}{\rho(y)}\dy \ \ge\ \left(\int_D g(y)\dy\right)^2,
\end{equation*}
where equality holds if and only if \(\rho(y)\propto g(y)\). Thus the variance-minimizing importance density is given by
\begin{equation*}
\rho^{*}(y)\ \propto\ g(y)\ =\ |\mathcal{R}(y)|^2.
\end{equation*}

The variance comparison depends on the second moment $\int_D g^2/\rho$. Normalizing $g^2$ therefore gives a useful reference density for this comparison. The following theorem provides a sufficient variance-reduction criterion for normalized densities that lie pointwise between this reference density and the uniform density:
\begin{theorem}
\label{thm:statistic_error}
Assume $0<\int_Dg^2(y)\dy<\infty$ and let \(\displaystyle\tilde{\rho}(y) := \frac{g^2(y)}{\int_D g^2(y)\dy}\). Suppose \(\rho(y)>0\) almost everywhere is a normalized probability density on $D$ such that for some measurable function \(\tau(y) \in [0,1]\),
\begin{equation}\label{eq4.1}
\rho(y) := (1 - \tau(y)) + \tau(y) \cdot \tilde{\rho}(y), \quad \text{a.e. on } D.
\end{equation}
Then the following inequality holds:
\begin{equation}\label{eq:mixture_upperbound}
	\int_D \frac{g^2(y)}{\rho(y)} \dy \le \int_D g^2(y)\dy,
\end{equation}
and consequently,
\begin{equation*}
	\mathrm{Var}_{\rho}\left(\frac{g(y)}{\rho(y)}\right) \le \mathrm{Var}_{\rho_u}\big(g(y)\big).
\end{equation*}
Moreover, the inequality is strict if the set where $0<\tau(y)<1$ and $\tilde\rho(y)\ne1$ has positive measure.
\end{theorem}

\begin{proof}
Writing $q=\tilde\rho$, the pointwise identity
\[
    (1-\tau)q+\tau-\frac{q}{\rho}
    =\frac{\tau(1-\tau)(q-1)^2}{\rho}\ge0
\]
gives
\begin{equation}\label{eq4.2}
    \frac{\tilde\rho}{\rho}
    \le (1-\tau)\tilde\rho+\tau.
\end{equation}
Integrating \eqref{eq4.1} over \(D\) and using the normalization of density function, we have
\begin{equation*}
\int_D \rho\dy = \int_D \left( (1 - \tau)  + \tau \tilde{\rho} \right)\dy = 1 + \int_D \tau(\tilde{\rho} - 1)\dy = 1,
\end{equation*}
which yields \(\int_D \tau(\tilde{\rho} - 1)\dy = 0\). Therefore, using \eqref{eq4.2}, we have
\begin{equation}\label{eq4.3}
\int_D \frac{\tilde{\rho}}{\rho}\dy \le \int_D \left( (1 - \tau)\tilde{\rho} + \tau \right) \dy =\int_D \tilde{\rho}\dy - \int_D \tau(\tilde{\rho} - 1)\dy= 1.
\end{equation}
By the definition of $\tilde{\rho}$, we have \(\displaystyle\int_D \frac{g^2}{\rho}\dy = \int_D g^2\dy\!\int_D \frac{\tilde{\rho}}{\rho}\dy\). This identity, together with \eqref{eq4.3}, yields \eqref{eq:mixture_upperbound}. The pointwise identity also gives strictness under the stated positive-measure condition. 
\end{proof}

In our framework, uniform computational samples induce the spatial density $\rho_t$ in \eqref{eq:coordinate_density}, as derived in Section~\ref{sec:lifting_variable}. With uniform time sampling, the corresponding density on $D$ is $\rho(t,\mathbf{x})=\rho_t(\mathbf{x})/T$.
In one dimension, equidistribution gives $\rho_t(x)=M(t,x)/\Theta(t)$. The monitor becomes large near sharp transition layers, concentrating physical samples where the solution is hardest to resolve. When the squared residual $g(y)=|\mathcal R(y)|^2$ is concentrated in these same layers, the geometric density aligns with the regions important for residual integration. Theorem~\ref{thm:statistic_error} specifies a sufficient condition for this alignment to reduce estimator variance. This provides an importance-sampling interpretation of the coordinate construction without explicit residual tracking.

\subsection{Gradient flow perspective: lifting accelerates training}
\label{sec:sec43}

We study the empirical training loss through the NTK and gradient-flow formulation for wide networks~\cite{jacot2018ntk,wang2022ntk}. Consider a linear flux $f(u)$, so that the PDE residual operator $\mathcal{R}$ is linear, and fix the lifting coordinate and collocation set. Let $y:=(t,\mathbf{x})\in D=[0,T]\times\Omega_P$. The training set consists of initial/boundary points $\{y_i^{ib}\}_{i=1}^{N_{ib}}$ with targets $\{u_i^{ib}\}$ and interior collocation points $\{y_i^{r}\}_{i=1}^{N_{r}}$ used to enforce the PDE residual.

The total empirical loss for the neural network $\tilde{u}(\cdot; \theta)$ is composed of a residual term and an initial-boundary term: $\ell(\theta) = \ell_r(\theta) + \ell_{ib}(\theta)$, where
\[
\ell_r(\theta)=\frac12\sum_{i=1}^{N_r}a_i\big(\mathcal{R}(\tilde{u})(y_i^r;\theta)\big)^2,\quad
\ell_{ib}(\theta)=\frac12\sum_{i=1}^{N_{ib}}b_i\big(\tilde{u}(y_i^{ib};\theta)-u_i^{ib}\big)^2.
\]
Here $a_i,b_i>0$ include the loss coefficients, sample normalization, and residual measure correction. Let \(s\ge 0\) represent the continuous training time. The parameter vector ${\theta}(s)$ evolves according to the gradient flow:
\[
\dot{\theta}(s)=-\nabla_\theta \ell(\theta(s)),\quad \theta(0)=\theta_0.
\]
We define the initial-boundary residual vector $r_{ib}(s) \in \mathbb{R}^{N_{ib}}$ and the PDE residual vector $r_r(s) \in \mathbb{R}^{N_r}$ by
\[
r_{ib}(s)
:=\left[\sqrt{b_i}\big(\tilde u(y_i^{ib};\theta(s))-u_i^{ib}\big)\right]_{i=1}^{N_{ib}},
\quad
r_r(s)
:=\left[\sqrt{a_i}\mathcal R(\tilde u)(y_i^r;\theta(s))\right]_{i=1}^{N_r}.
\]
The full residual vector is then given by $r(s) = [r_{ib}(s)^\top, r_r(s)^\top]^\top \in \mathbb{R}^{N}$, where $N = N_{ib} + N_r$.

To analyze the training dynamics, we define the Jacobian $J(s) = \partial r / \partial \theta |_{\theta(s)}$ in block form, corresponding to the initial-boundary and PDE residual terms:
\[
J(s) = \begin{bmatrix} J_{ib}(s) \\ J_r(s) \end{bmatrix}, \quad
\begin{aligned}
	\text{where} \quad [J_{ib}(s)]_{i,:} &= \sqrt{b_i}\,\nabla_\theta \tilde{u}(y_i^{ib}; \theta(s))^\top, \\
	[J_r(s)]_{i,:} &= \sqrt{a_i}\,\nabla_\theta \mathcal{R}(\tilde{u})(y_i^r; \theta(s))^\top.
\end{aligned}
\]
Combining this with the gradient flow $\dot{\theta}(s) = -J(s)^\top r(s)$, the chain rule gives the residual dynamics:
\begin{equation}\label{eq:residual_dynamics}
	\dot{r}(s) = -K(s) r(s), \quad\text{where}\,
	K(s) := J(s) J(s)^\top=
	\begin{bmatrix}
		K_{ib,ib}(s) & K_{ib,r}(s) \\
		K_{r,ib}(s) & K_{r,r}(s)
	\end{bmatrix}.
\end{equation}
The matrix $K(s)$ is the weighted NTK, which governs the convergence behavior of the network. The chain-rule identity \eqref{eq:residual_dynamics} also holds for differentiable nonlinear residuals.

Under the infinite-width limit for single-hidden-layer networks with appropriate initialization, the NTK $K(s)$ remains approximately constant throughout training~\cite{jacot2018ntk, wang2022ntk}, i.e.,
\[
K(s)\approx K(0)=:K^*=
	\begin{bmatrix}
		K_{ib,ib}^* & K_{ib,r}^* \\
		K_{r,ib}^* & K_{r,r}^*
	\end{bmatrix}.
\]

Consequently, the residual dynamics simplify to $
r(s)\approx \exp(-K^* s)\, r(0)$.
Now, let \(J^* = J(0)\) have the singular value decomposition \(J^* = U \Sigma V^\top\), where \(\Sigma = \operatorname{diag}(\sigma_1, \dots, \sigma_m)\) and \(m = \operatorname{rank}(J^*)\). Then,
\[
K^* = J^* J^{*\top} = U \Sigma^2 U^\top,
\]
which implies that in the coordinate system defined by \(U\), the residuals evolve as:
\[
U^\top r(s) \approx \exp(-\Sigma^2 s)\, U^\top r(0).
\]
For the compact SVD, the remaining component satisfies $(I-UU^\top)r(s)\approx(I-UU^\top)r(0)$ under the frozen kernel. Thus each positive singular value $\sigma_i$ governs the exponential decay rate of the residual component along its corresponding singular direction~\cite{wang2022ntk}.

Heuristically, adaptive sampling concentrates the collocation points \(\{y_i^r\}\) in regions where the residual is large, which can increase the correlation among the corresponding tangent features (i.e., the rows of \(J^*\)). When this clustering compresses the tail of the singular-value spectrum, smaller singular values slow the decay of the associated residual modes. Therefore, although adaptive sampling enhances local approximation near sharp features, it can also impede the overall convergence rate.

To mitigate this effect, we introduce a lifting transformation that augments each collocation point with an auxiliary coordinate, forming an extended set \(\{(y_i^r, \mathbf{z}_i^r)\}\). By expanding the sampling space, the lifting operation separates previously clustered points in the lifted domain, providing additional coordinate variation that helps reduce correlations among the tangent features in \(J^*\). This explains how lifting can improve the lower part of the NTK spectrum and accelerate the decay of the corresponding residual modes. As a result, the lifting strategy preserves the advantage of sample concentration in high-residual regions while alleviating the adverse effect on training dynamics.

\section{Numerical experiments}
\label{sec:numerical_experiments}

\setcounter{topnumber}{3}
\setcounter{bottomnumber}{2}
\setcounter{totalnumber}{5}
\renewcommand{\topfraction}{0.92}
\renewcommand{\bottomfraction}{0.85}
\renewcommand{\textfraction}{0.06}
\renewcommand{\floatpagefraction}{0.65}

In this section, we evaluate the approximation accuracy, training convergence, and computational cost of TAL-PINN for hyperbolic conservation laws near the inviscid limit.

After the common experimental settings, a mechanism study uses coordinates constructed from the reference solution to distinguish the effects of adaptive sampling and lifting through sample distributions, NTK spectra, and loss decay. Stationary and moving Burgers tests then assess the autonomous method with coordinates constructed from coarse solutions; an ablation study in the stationary case further separates the contributions of lifting, adaptive sampling, and viscosity continuation. The one-dimensional Euler tests examine the extension from scalar equations to systems of conservation laws, while the two-dimensional Burgers test evaluates the method in two spatial dimensions. Finally, we compare the time spent on solution training, mesh generation, and coordinate fitting to assess the computational cost of the accuracy gains.

\subsection{Experimental setup}
\label{sec:exp_setup}

All experiments are implemented in JAX and jaxlib (version 0.9.2) and run on an NVIDIA Tesla V100-PCIE GPU with 16 GB memory. Table~\ref{tab:exp_architecture} lists the architectures of the solution and coordinate networks. Each network begins with a trainable Fourier feature layer~\cite{wang2021fourier},
\[
    \Phi_F(\mathbf{v})=\sin\!\left(2\pi(W_F\mathbf{v}+b_F)\right),
\]
followed by fully connected layers, with SiLU activations in the nonlinear hidden layers. The Fourier parameters $W_F$ and $b_F$ are optimized together with the remaining network parameters. The TAL-PINN solution network takes $(t,\mathbf{x},\mathbf{z}(t,\mathbf{x}))$ as input, whereas the coordinate network and the baseline solution networks take $(t,\mathbf{x})$.

\begin{table}[!htb]
\centering
\caption{Network architectures. $d_{\mathrm{in}}$ and $d_{\mathrm{out}}$ are the input and output dimensions, and $m_F$ is the Fourier feature width. The hidden widths for Euler refer to its shared trunk; each of its three branches has two hidden layers of width 40 and one scalar output. Parameter counts include the trainable Fourier layer. The Burgers baselines include the ablation variants.}
\label{tab:exp_architecture}
\small
\setlength{\tabcolsep}{3pt}
\begin{tabular}{@{}llrrcrr@{}}
\toprule
Problem & Network & $d_{\mathrm{in}}$ & $m_F$ & Hidden & $d_{\mathrm{out}}$ & Parameters\\
\midrule
1D Burgers & TAL solution & 3 & 20 & $40$--$40$ & 1 & 2601\\
 & Baseline solution & 2 & 20 & $40$--$40$ & 1 & 2581\\
 & Coordinate & 2 & 20 & $40$--$40$ & 1 & 2581\\
\midrule
1D Euler & TAL solution & 3 & 20 & $40$--$40$ & 3 & 12520\\
 & Vanilla solution & 2 & 40 & $40$--$40$ & 3 & 13360\\
 & Coordinate & 2 & 20 & $40$--$40$ & 1 & 2581\\
\midrule
2D Burgers & TAL solution & 5 & 40 & $40$--$40$ & 1 & 3561\\
 & Vanilla solution & 3 & 40 & $40$--$40$ & 1 & 3481\\
 & Coordinate & 3 & 40 & $40$--$40$ & 2 & 3522\\
\bottomrule
\end{tabular}
\end{table}

All networks are trained using SOAP~\cite{vyas2024soap}, with learning rate $3\times10^{-3}$, momentum coefficients $(0.95,0.95)$, weight decay $0.01$, and a preconditioner update every 10 iterations. The learning rate is fixed within each training stage. At each viscosity reduction, the solution-network parameters are carried over and the optimizer state is restarted.

In the mechanism study, the coordinates are constructed from the reference solution at a fixed viscosity. In the autonomous TAL-PINN experiments, adaptive coordinates are constructed from the preceding coarse solution according to Algorithm~\ref{alg:Lifting_framework} and its multistage extension. Each coordinate fit minimizes \eqref{eq:loss_coordnet} on the mesh data $\{(t_i,\mathbf{x}_i,\boldsymbol{\xi}_i)\}$ for $10^4$ iterations with batches of $10^4$ nodal samples. The fitted parameters $\theta_\xi^*$ remain fixed during the subsequent optimization of the solution-network parameters $\theta$.

Residual, initial, and boundary collocation points are resampled at each solution-training iteration, with the batch sizes and loss weights given in Table~\ref{tab:exp_batches}. For adaptive sampling, the residual loss uses self-normalized measure-correction weights. The ablation variants apply viscosity continuation to a vanilla PINN with uniform or adaptive sampling. Table~\ref{tab:exp_schedule} specifies the viscosities and solution-training budgets at successive stages, together with the number of coordinate fits.

\begin{table}[!htb]
\centering
\caption{Collocation batches and loss weights. $N_r$ and $N_{\mathrm{ic}}$ count residual and initial points per iteration. $N_{\mathrm{bc}}^{\mathrm{side}}$ counts points at each endpoint for 1D Burgers and periodic pairs in each spatial direction for 2D Burgers, giving $N_{\mathrm{bc}}=2N_{\mathrm{bc}}^{\mathrm{side}}$ boundary constraints in either case. Euler uses the PDE and initial losses.}
\label{tab:exp_batches}
\small
\setlength{\tabcolsep}{5pt}
\begin{tabular}{@{}lrrrl@{}}
\toprule
Problem & $N_r$ & $N_{\mathrm{ic}}$ & $N_{\mathrm{bc}}^{\mathrm{side}}$ & $(w_r,w_{\mathrm{ic}},w_{\mathrm{bc}})$\\
\midrule
1D Burgers & 10000 & 1000 & 1000 & $(1,10,10)$\\
1D Euler & 10000 & 1000 & --- & $(1,100,0)$\\
2D Burgers & 50000 & 5000 & 5000 & $(1,10,10)$\\
\bottomrule
\end{tabular}
\end{table}

\begin{table}[!htb]
\centering
\caption{Training schedules. Entries separated by slashes refer to successive stages. PDE iterations count solution-network updates and are given in thousands. The last column counts coordinate-network fits.}
\label{tab:exp_schedule}
\small
\setlength{\tabcolsep}{3pt}
\begin{tabular}{@{}llccc@{}}
\toprule
Problem & Method & Viscosities & PDE iterations & Coord. fits\\
\midrule
Mechanism & Three methods & $10^{-4}$ & $10$ & 1 (TAL)\\
\midrule
1D Burgers & TAL-PINN & $10^{-2}/10^{-3}/10^{-4}$ & $2/2/4$ & 2\\
(both shocks) & Vanilla PINN & $10^{-3}$ & $8$ & 0\\
 & Vanilla PINN & $10^{-4}$ & $8$ & 0\\
 Stationary & Ablation variants & $10^{-2}/10^{-3}/10^{-4}$ & $2/6/4$ & 0\\
\midrule
Sod & TAL-PINN & $10^{-3}/10^{-4}$ & $10/10$ & 1\\
 & Vanilla PINN & $10^{-4}$ & $20$ & 0\\
Lax & TAL-PINN & $3\!\times\!10^{-3}/3\!\times\!10^{-4}$ & $10/10$ & 1\\
 & Vanilla PINN & $3\times10^{-4}$ & $20$ & 0\\
\midrule
2D Burgers & TAL-PINN & $10^{-2}/10^{-3}$ & $4/10$ & 1\\
 & Vanilla PINN & $10^{-3}$ & $14$ & 0\\
\bottomrule
\end{tabular}
\end{table}

\Needspace{8\baselineskip}
Since the objective is to approximate inviscid solutions, all test errors are measured against an inviscid reference $u^*$. For a test set $\mathcal{T}$, we report the relative $L^2$ error,
\begin{equation}
\operatorname{RL2}(u) = \left( \frac{\sum_{i\in\mathcal{T}}|u_i-u_i^*|^2} {\sum_{i\in\mathcal{T}}|u_i^*|^2} \right)^{1/2}.
\label{eq:exp_metrics}
\end{equation}
For Euler, the relative $L^2$ error is evaluated separately for density, velocity, and pressure. All experiments are repeated ten times, and the result tables and convergence curves report mean values.

\FloatBarrier

\subsection{Mechanism verification}
\label{sec:exp_mechanism}

We consider the one-dimensional inviscid Burgers equation
\begin{equation}
    u_t+uu_x=0,\qquad (t,x)\in[0,1]\times[0,1],
    \label{eq:burgers}
\end{equation}
with $u(0,x)=\sin(2\pi x)$ and homogeneous Dirichlet conditions $u(t,0)=u(t,1)=0$. Although the initial condition is smooth, a stationary shock develops at $x=0.5$ when $t=(2\pi)^{-1}$. In the PINN formulation, we solve the regularized equation
\begin{equation}
    u_t+uu_x=\nu u_{xx},
    \label{eq:burgers_viscous}
\end{equation}
with $\nu=10^{-4}$. The inviscid reference is computed using a conservative Godunov finite-volume scheme~\cite{leveque2002finite}, with test data on a $513\times513$ space--time grid for the 1D Burgers experiments.

We first illustrate the coordinate-induced sampling mechanism discussed in Section~\ref{sec:sec42}. Using spatial gradients computed from the discrete reference solution, we evaluate the monitor function
\[
    M(t,x)=\sqrt{1+(u_x^*(t,x))^2}.
\]
After smoothing the monitor function in space, we construct the adaptive coordinates using the one-dimensional equidistribution formula \eqref{eq:equidistribution_1d_map}. Uniform samples in $\Omega_C$ are mapped to physical residual points through the interpolated forward map $x_h(t,\xi)$. In TAL-PINN, the fitted coordinate $\tilde{\xi}(t,x;\theta_\xi^*)$ also serves as the lifting input.

\begin{figure}[!htb]
\centering
\begin{subfigure}[t]{0.480\linewidth}
    \centering
    \includegraphics[width=\linewidth]{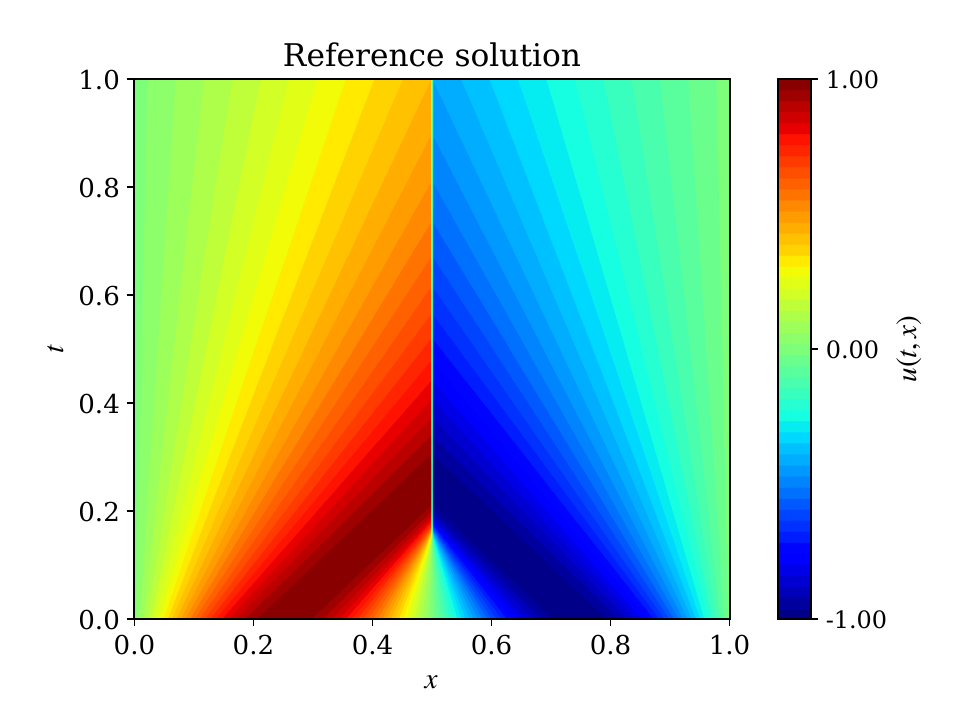}
    \caption{Inviscid reference}
\end{subfigure}
\hfill
\begin{subfigure}[t]{0.480\linewidth}
    \centering
    \includegraphics[width=\linewidth]{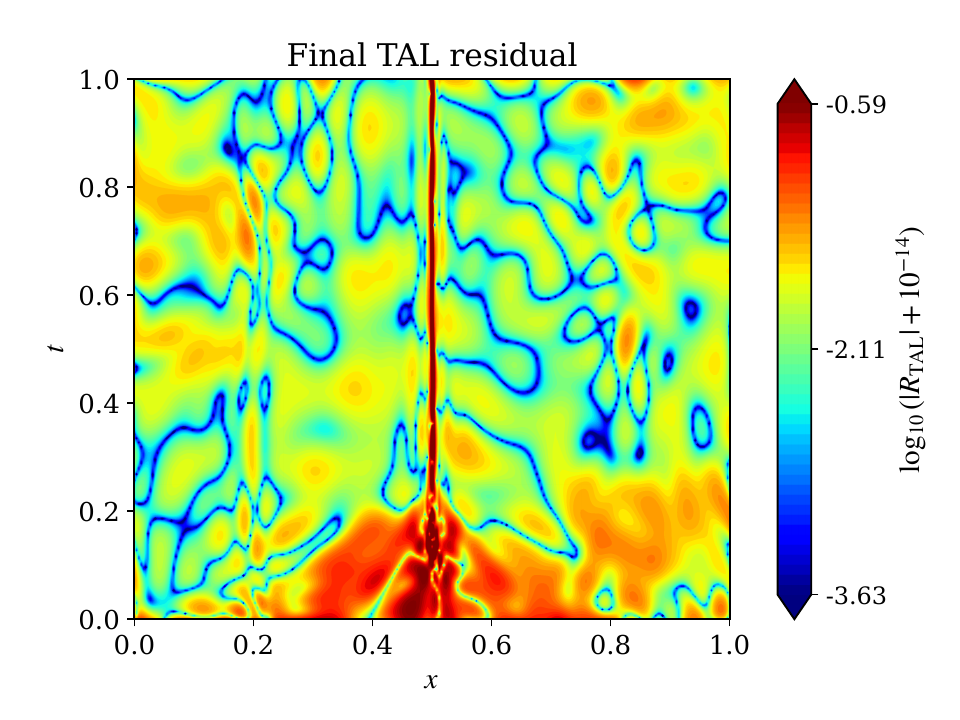}
    \caption{Final residual}
\end{subfigure}
\par\medskip
\begin{subfigure}[t]{0.480\linewidth}
    \centering
    \includegraphics[width=\linewidth]{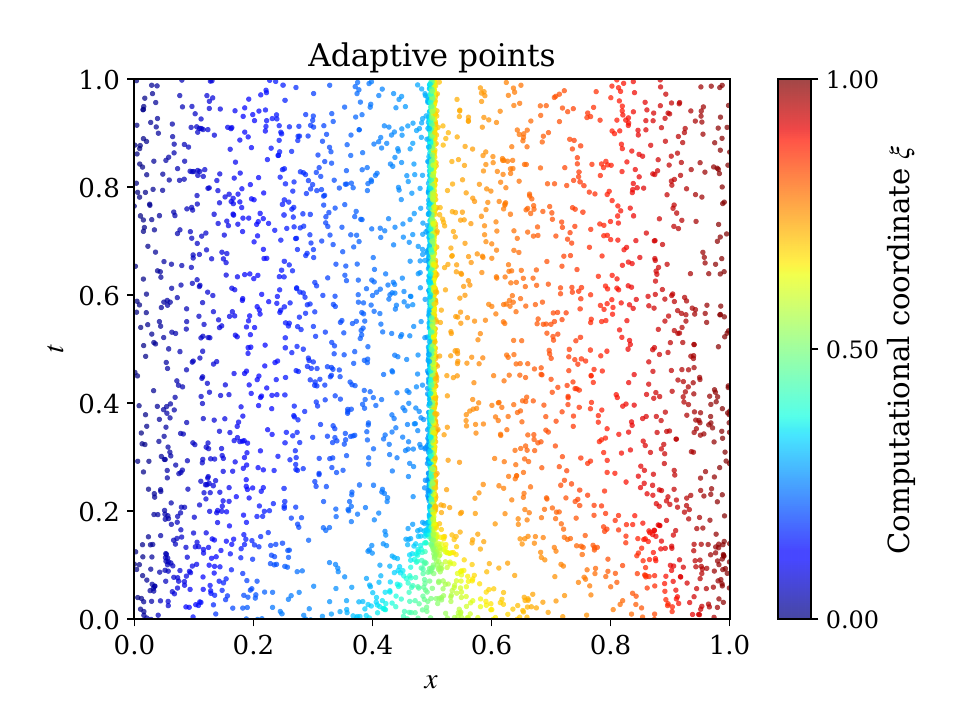}
    \caption{Samples in physical space}
\end{subfigure}
\hfill
\begin{subfigure}[t]{0.480\linewidth}
    \centering
    \includegraphics[width=\linewidth]{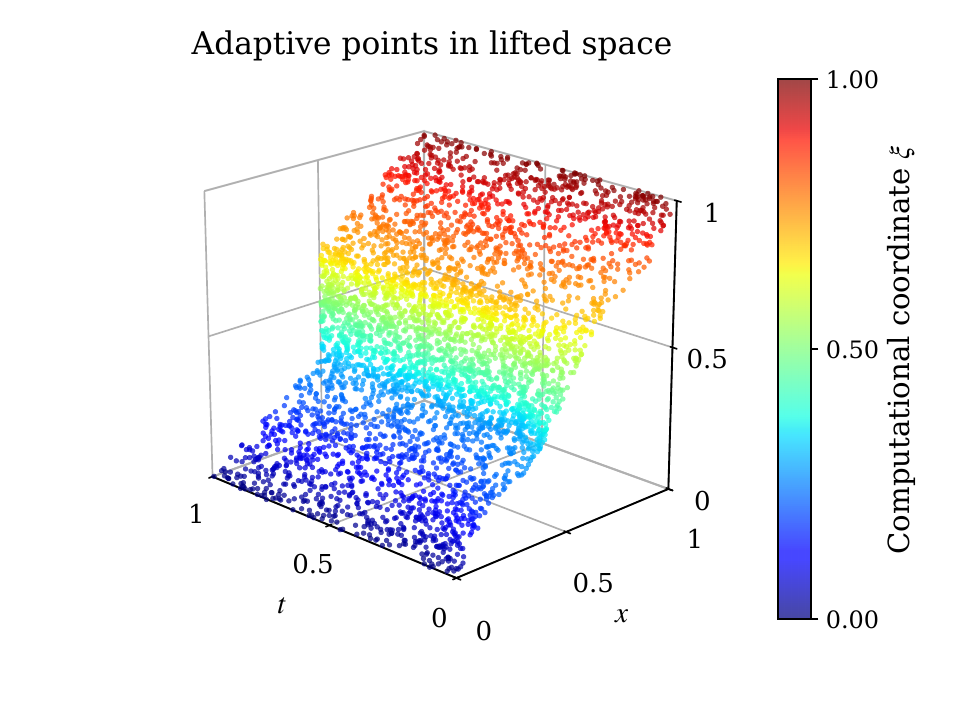}
    \caption{Samples in lifted space}
\end{subfigure}
\caption{Adaptive sampling and lifting for the stationary Burgers shock at $\nu=10^{-4}$: (a) inviscid reference solution; (b) final TAL-PINN residual magnitude on a logarithmic scale; (c) adaptive samples in physical space; and (d) the corresponding samples in lifted space. Colors in (c,d) indicate the computational coordinate $\xi$. Results are shown for one run.}
\label{fig:exp_mechanism_sampling}
\end{figure}

Figure~\ref{fig:exp_mechanism_sampling}(a,b) shows the reference solution and the final TAL-PINN residual. Large residuals occur near the stationary shock and its formation region. Panel (c) shows that uniform sampling in computational coordinates produces a higher density of physical points near these regions. This spatial agreement between sample concentration and residual localization supports the importance-sampling mechanism in Section~\ref{sec:sec42}: the adaptive coordinate directs training points toward regions that contribute strongly to the residual integral.

Panel (d) displays the corresponding samples in the lifted space $(t,x,\xi)$. Near the shock, the auxiliary coordinate varies rapidly over a narrow physical region, separating nearby physical points along the additional coordinate direction. This view illustrates how the adaptive coordinate encodes the transition geometry in the lifted representation, complementing its role in adaptive sampling.

We next examine the optimization mechanism in Section~\ref{sec:sec43} by comparing a vanilla PINN with uniform sampling, an adaptive-sampling PINN without lifting, and TAL-PINN. The two adaptive methods use the same interpolated forward map $x_h(t,\xi)$ to generate physical residual points. TAL-PINN additionally uses $\tilde{\xi}(t,x;\theta_\xi^*)$ as its lifting input, allowing the effect of lifting to be assessed under a common adaptive sampling distribution.

We compare the total NTK spectrum at initialization and the subsequent loss decay. The NTK is computed with respect to the solution-network parameters, keeping the fitted coordinate fixed for TAL-PINN. We use 1000 residual points, 100 initial points, and 200 boundary points split equally between the two endpoints. The kernel incorporates the loss weights, sample normalization, and adaptive-sampling correction, following Section~\ref{sec:sec43}.

\begin{figure}[!htb]
\centering
\begin{subfigure}[t]{0.480\linewidth}
    \centering
    \includegraphics[width=\linewidth]{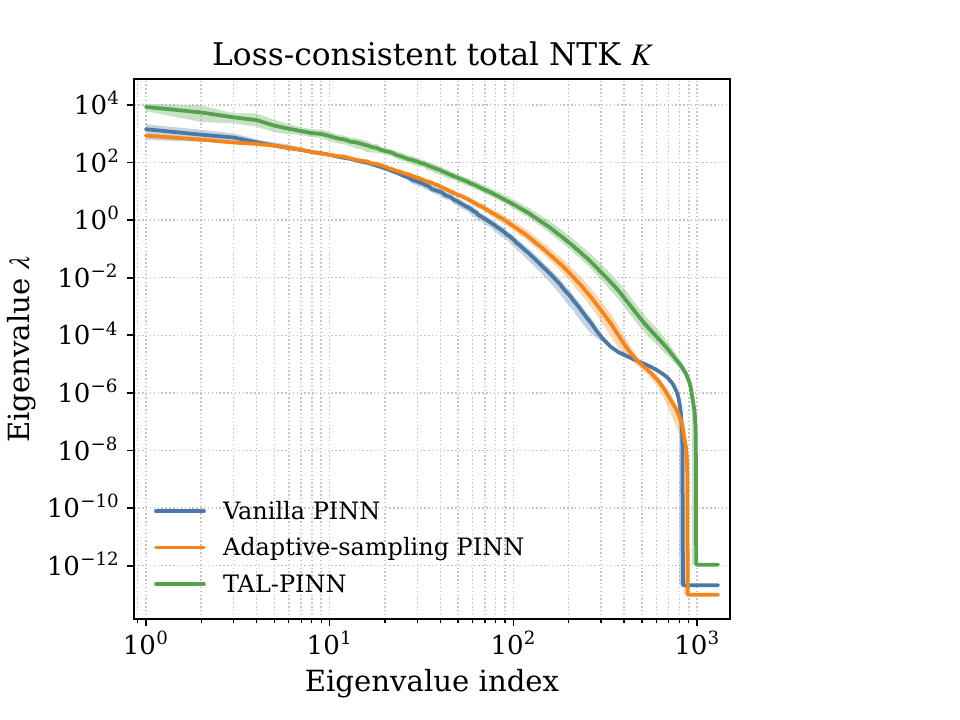}
    \caption{Total NTK}
\end{subfigure}
\hfill
\begin{subfigure}[t]{0.480\linewidth}
    \centering
    \includegraphics[width=\linewidth]{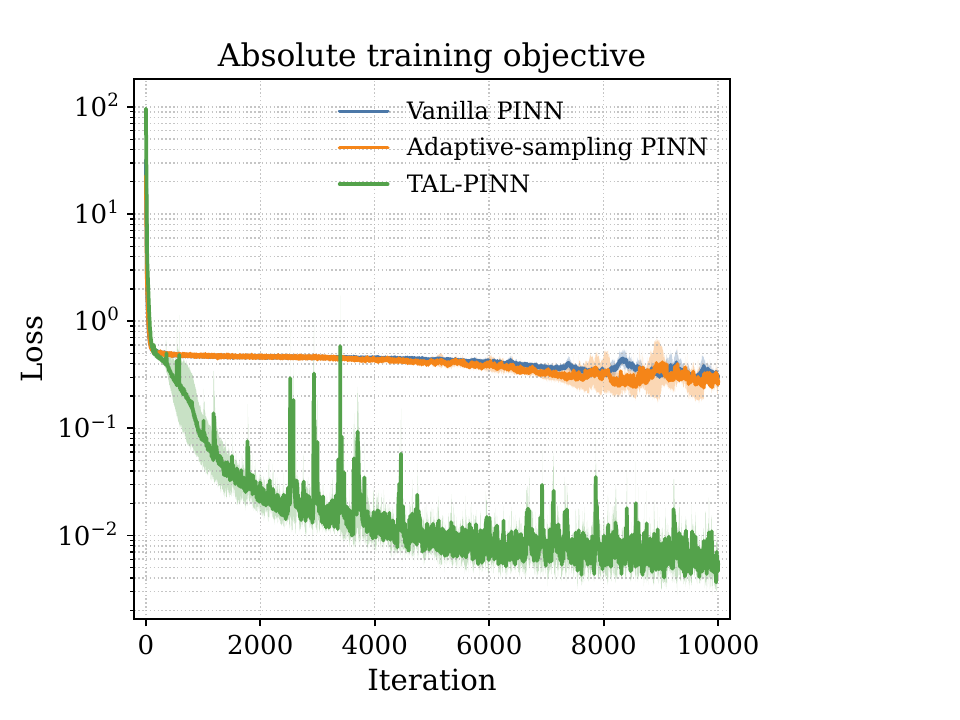}
    \caption{Training objective}
\end{subfigure}
\caption{Optimization mechanism at $\nu=10^{-4}$: (a) total NTK spectra at initialization and (b) training-objective histories for vanilla PINN, adaptive-sampling PINN, and TAL-PINN. Shading indicates the minimum--maximum range.}
\label{fig:exp_mechanism_optimization}
\end{figure}

Figure~\ref{fig:exp_mechanism_optimization}(a) shows that TAL-PINN has substantially larger eigenvalues over a broad portion of the total NTK spectrum, consistent with faster residual-mode decay in the gradient-flow analysis. This spectral improvement is reflected in the training objectives in panel (b): TAL-PINN rapidly reduces the loss, while the vanilla and adaptive-sampling PINNs plateau at similar levels. At this viscosity, adaptive sampling alone therefore gives little improvement in optimization, despite concentrating points near the shock. With the same adaptive sampling distribution, using the coordinate as an additional network input improves the representation of the steep solution and the resulting optimization behavior.

\subsection{1D Burgers equation with a stationary shock}
\label{sec:burgers_1d}

We now apply the autonomous multistage procedure to the same stationary-shock problem. TAL-PINN starts from the identity coordinate and trains successively at $\nu=10^{-2}$, $10^{-3}$, and $10^{-4}$. After each of the first two stages, the current solution defines the monitor $M=\sqrt{1+u_x^2}$, from which the mesh and the next lifting coordinate are constructed. No reference solution is used in these updates. The coordinate loss \eqref{eq:loss_coordnet} uses $w_J=0.01$ and $\delta_J=0.05$ to discourage loss of monotonicity.

\begin{figure}[!htb]
\centering
\begin{subfigure}[t]{0.480\linewidth}
    \centering
    \includegraphics[width=\linewidth]{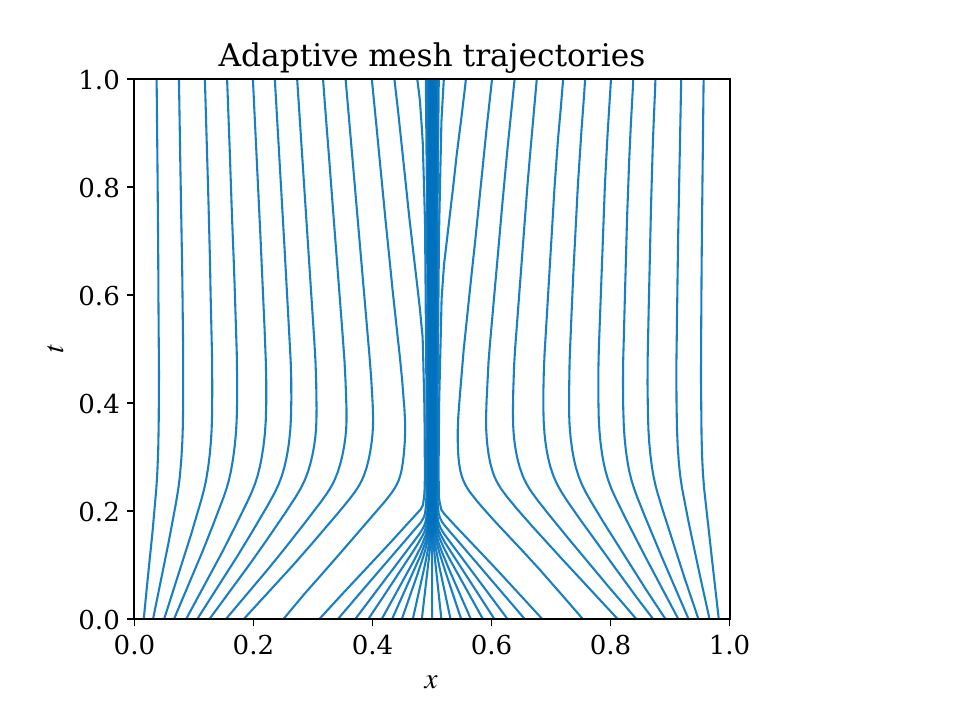}
    \caption{Adaptive mesh}
\end{subfigure}
\hfill
\begin{subfigure}[t]{0.480\linewidth}
    \centering
    \includegraphics[width=\linewidth]{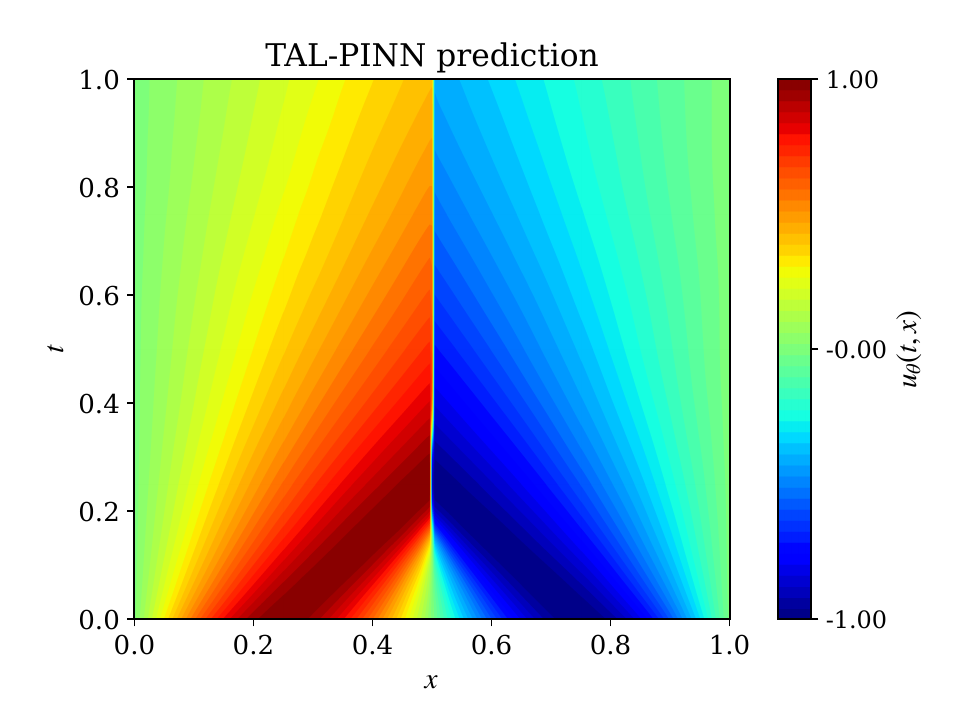}
    \caption{TAL-PINN prediction}
\end{subfigure}
\par\medskip
\begin{subfigure}[t]{0.480\linewidth}
    \centering
    \includegraphics[width=\linewidth]{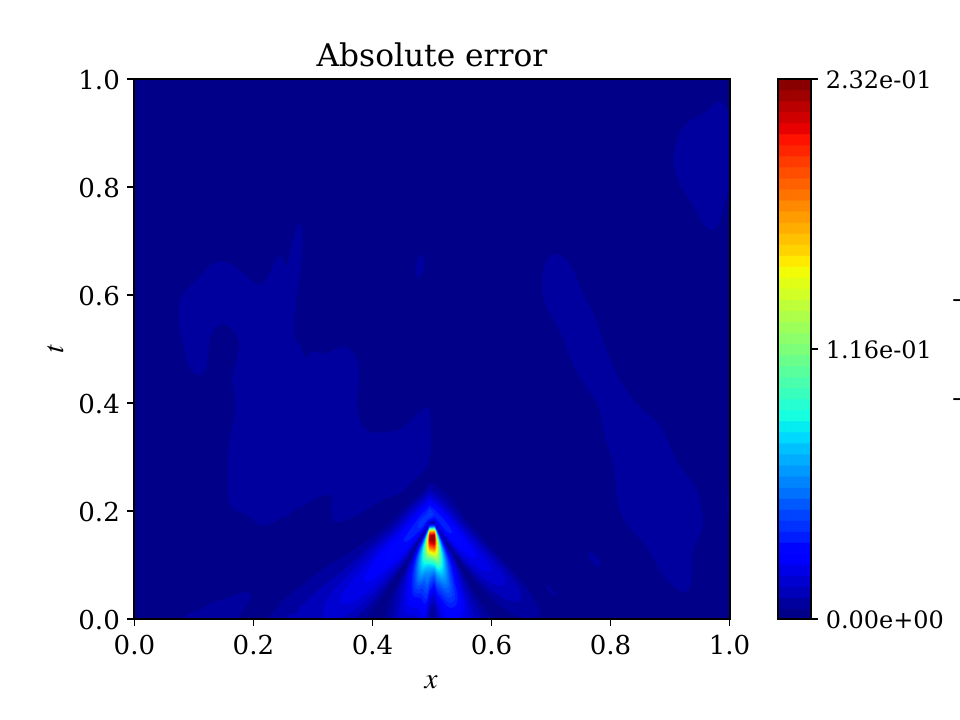}
    \caption{Absolute error}
\end{subfigure}
\hfill
\begin{subfigure}[t]{0.480\linewidth}
    \centering
    \includegraphics[width=\linewidth]{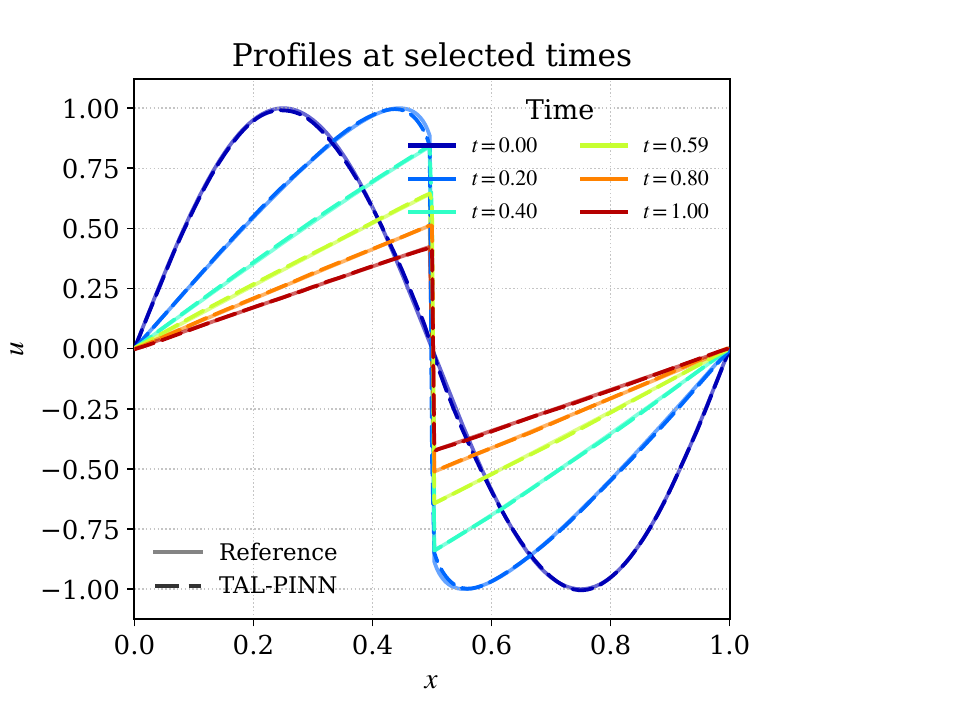}
    \caption{Time profiles}
\end{subfigure}
\caption{Stationary Burgers shock: adaptive mesh, TAL-PINN prediction, absolute error, and profiles at selected times for one trained model.}
\label{fig:exp_static_solution}
\end{figure}

For comparison, vanilla PINNs are trained at fixed viscosities $10^{-3}$ and $10^{-4}$. All three methods use 8000 PDE iterations. Table~\ref{tab:exp_static} shows that TAL-PINN reduces the relative $L^2$ error at $\nu=10^{-4}$ from $4.48\times10^{-1}$ to $1.69\times10^{-2}$, while also improving upon the larger-viscosity baseline. The prediction and error in Figure~\ref{fig:exp_static_solution}(b,c) show why this gain matters: TAL-PINN captures the steep shock while preserving the smooth solution on either side. The remaining error is localized near the transition rather than spread throughout the domain. The larger vanilla-PINN error at $\nu=10^{-4}$ is consistent with the error decomposition in Section~\ref{sec:theoretical}. A smaller viscosity sharpens the shock layer and can make neural approximation and optimization more difficult. Moreover, the residual-to-solution bound in Lemma~\ref{thm:apost_L2} deteriorates as $\nu$ decreases. Thus, under a fixed training budget, the benefit of reduced viscous approximation error can be outweighed by increased neural approximation error.

\begin{table}[!htb]
\centering
\caption{Stationary Burgers shock: mean final relative $L^2$ errors.}
\label{tab:exp_static}
\small
\begin{tabular}{@{}lcc@{}}
\toprule
Method & Final $\nu$ & Relative $L^2$\\
\midrule
Vanilla PINN & $10^{-3}$ & $6.22\times10^{-2}$\\
Vanilla PINN & $10^{-4}$ & $4.48\times10^{-1}$\\
TAL-PINN & $10^{-4}$ & $1.69\times10^{-2}$\\
\bottomrule
\end{tabular}
\end{table}

\begin{figure}[!htb]
\centering
\begin{subfigure}[t]{0.480\linewidth}
    \centering
    \includegraphics[width=\linewidth]{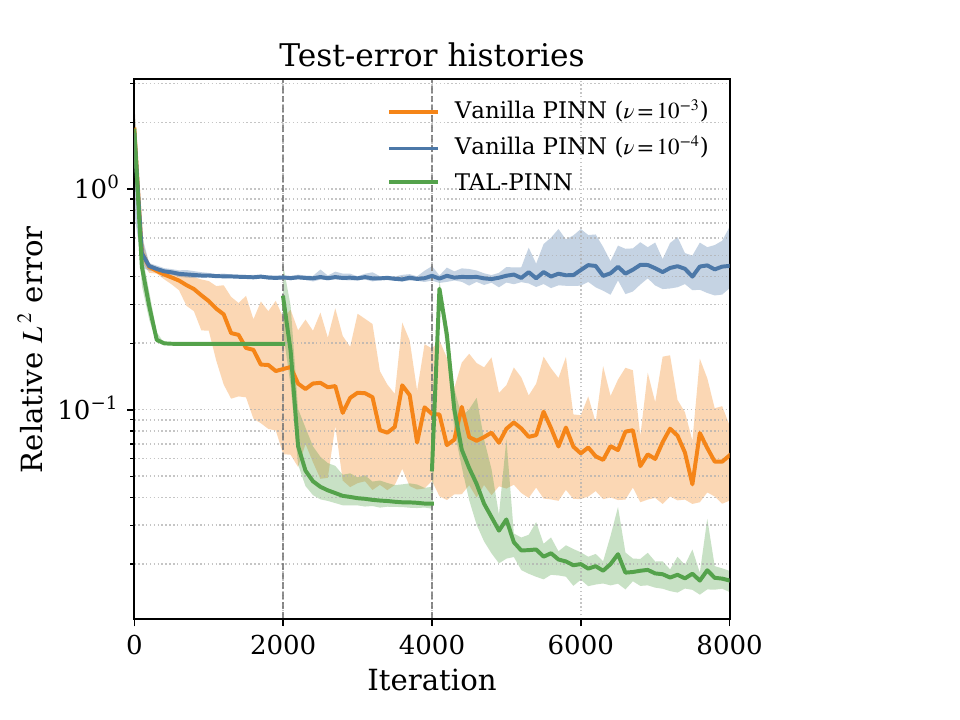}
    \caption{Comparison with vanilla PINNs}
\end{subfigure}
\hfill
\begin{subfigure}[t]{0.480\linewidth}
    \centering
    \includegraphics[width=\linewidth]{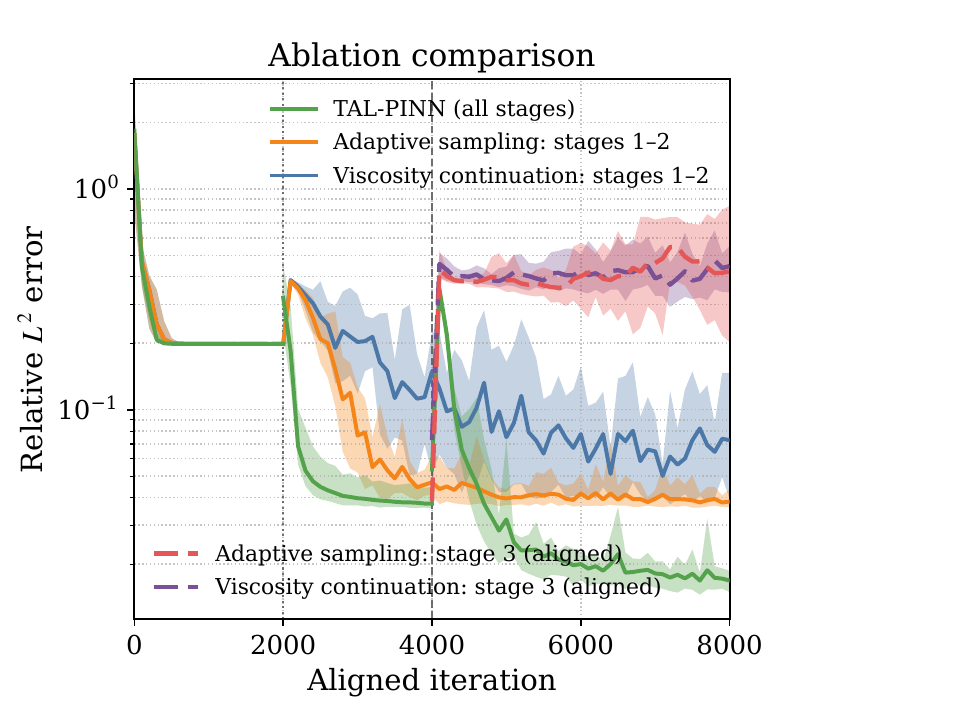}
    \caption{Ablation comparison}
\end{subfigure}
\caption{Test-error histories for the stationary Burgers shock: (a) comparison with fixed-viscosity vanilla PINNs; (b) comparison with the ablation variants. Curves show mean relative $L^2$ errors, with minimum--maximum shading. In (a), TAL-PINN changes viscosity after 2000 and 4000 PDE iterations. In (b), ``Adaptive sampling'' denotes AS+VC; the final stages of VC and AS+VC (actual iterations 8000--12000) are aligned with TAL-PINN's final stage on 4000--8000.}
\label{fig:exp_static_convergence}
\end{figure}

The adaptive mesh in Figure~\ref{fig:exp_static_solution}(a) concentrates around $x=0.5$, and the profiles in panel (d) closely follow the inviscid reference as the shock develops. Figure~\ref{fig:exp_static_convergence}(a) shows a sustained reduction in test error through the viscosity stages, whereas direct training at the smallest viscosity stalls at a much larger error. Thus the multistage lifting strategy enables the network to benefit from reduced viscosity and recover a sharper, more accurate solution.

We further assess the contributions of viscosity continuation (VC), adaptive sampling (AS), and lifting by comparing TAL-PINN with two vanilla-PINN baselines: VC with uniform sampling and AS+VC with adaptive sampling. AS+VC uses the same adaptive sampling procedure and measure-correction weights as TAL-PINN, while its solution network takes only $(t,x)$ as input. Both variants follow the same viscosity sequence as TAL-PINN, with a longer middle stage, giving 12000 PDE iterations in total compared with 8000 for TAL-PINN.

Table~\ref{tab:exp_ablation} and Figure~\ref{fig:exp_static_convergence}(b) show that AS+VC improves upon VC and reaches an accuracy comparable to TAL-PINN at $\nu=10^{-3}$, using twice as many cumulative PDE updates. When the viscosity decreases to $10^{-4}$, the errors of both ablation variants increase, whereas TAL-PINN continues to improve. These results show that lifting enables further accuracy gains as the viscosity is reduced, beyond those obtained with viscosity continuation and adaptive sampling alone.

\begin{table}[!htb]
\centering
\caption{Stationary-shock ablation: mean relative $L^2$ errors at the end of each indicated stage. Cumulative iterations count PDE updates.}
\label{tab:exp_ablation}
\small
\setlength{\tabcolsep}{4pt}
\begin{tabular}{@{}lcrc@{}}
\toprule
Method & $\nu$ & Iterations & Relative $L^2$\\
\midrule
VC & $10^{-3}$ & 8000 & $7.28\times10^{-2}$\\
AS+VC & $10^{-3}$ & 8000 & $3.84\times10^{-2}$\\
TAL-PINN & $10^{-3}$ & 4000 & $3.76\times10^{-2}$\\
\midrule
VC & $10^{-4}$ & 12000 & $4.47\times10^{-1}$\\
AS+VC & $10^{-4}$ & 12000 & $4.23\times10^{-1}$\\
TAL-PINN & $10^{-4}$ & 8000 & $1.69\times10^{-2}$\\
\bottomrule
\end{tabular}
\end{table}

\subsection{1D Burgers equation with a moving shock}
\label{sec:burgers_1d_moving}

To extend the evaluation to a moving shock, we revisit \eqref{eq:burgers} with the same homogeneous Dirichlet conditions and the initial condition
\[
    u(0,x)=\sin(2\pi x)+\frac12\sin(\pi x).
\]
The additional low-frequency component produces a curved shock trajectory in space--time. We use the same network architectures, monitor, collocation batches, and three-stage training schedule as in the stationary case. The adaptive coordinate is updated from the learned solution after the first and second stages, allowing both the lifting input and the sampling distribution to adapt to the evolving shock geometry.

\begin{figure}[!htbp]
\centering
\begin{subfigure}[t]{0.480\linewidth}
    \centering
    \includegraphics[width=\linewidth]{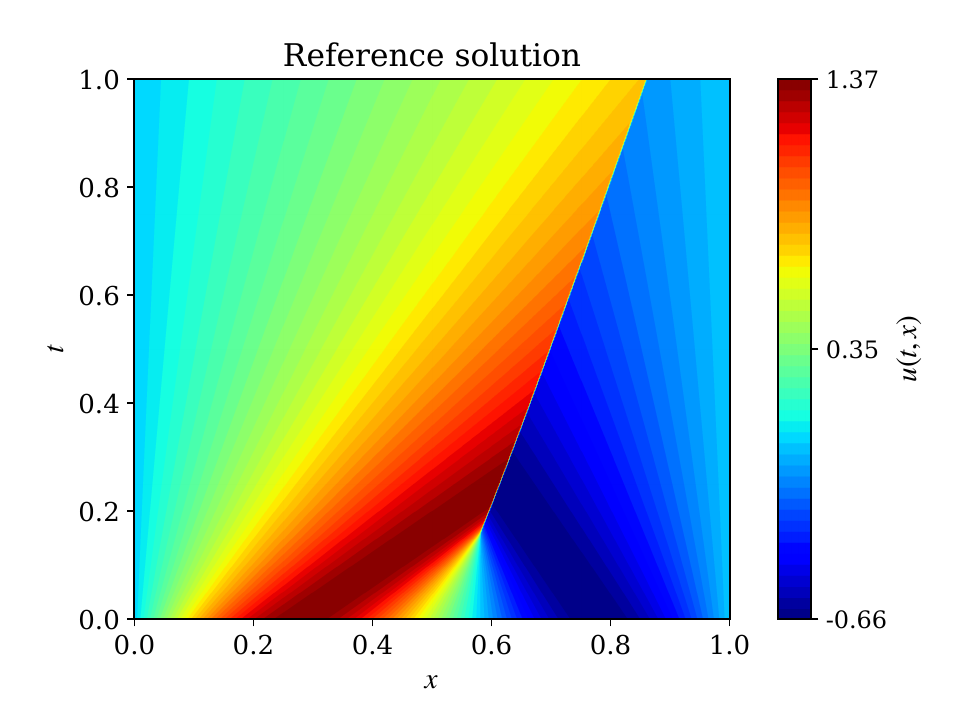}
    \caption{Inviscid reference}
\end{subfigure}\hfill
\begin{subfigure}[t]{0.480\linewidth}
    \centering
    \includegraphics[width=\linewidth]{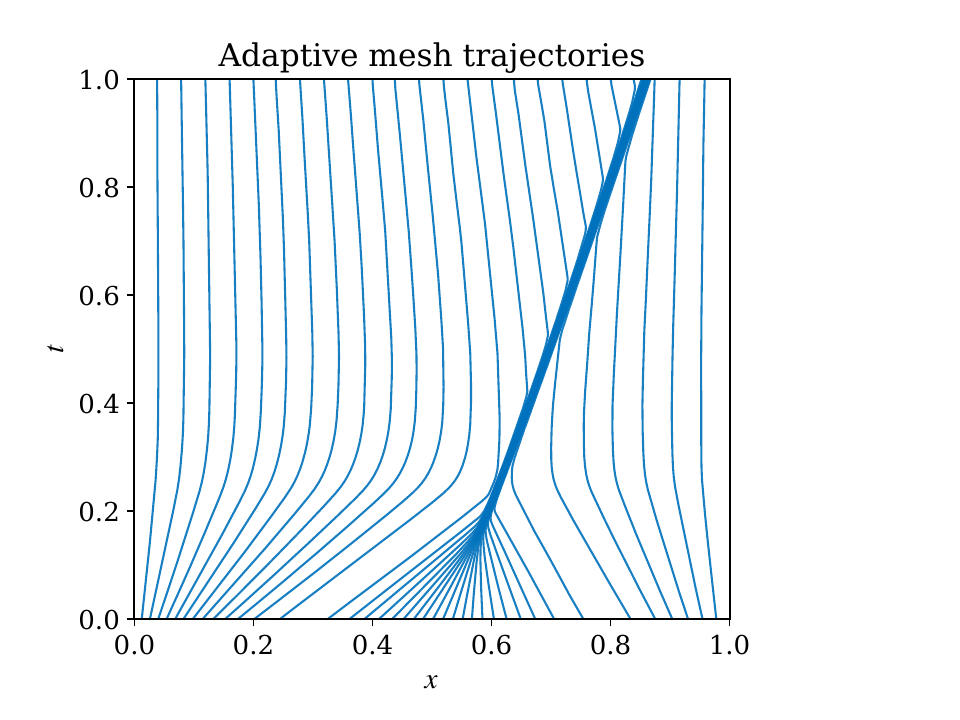}
    \caption{Adaptive mesh}
\end{subfigure}

\par\medskip

\begin{subfigure}[t]{0.480\linewidth}
    \centering
    \includegraphics[width=\linewidth]{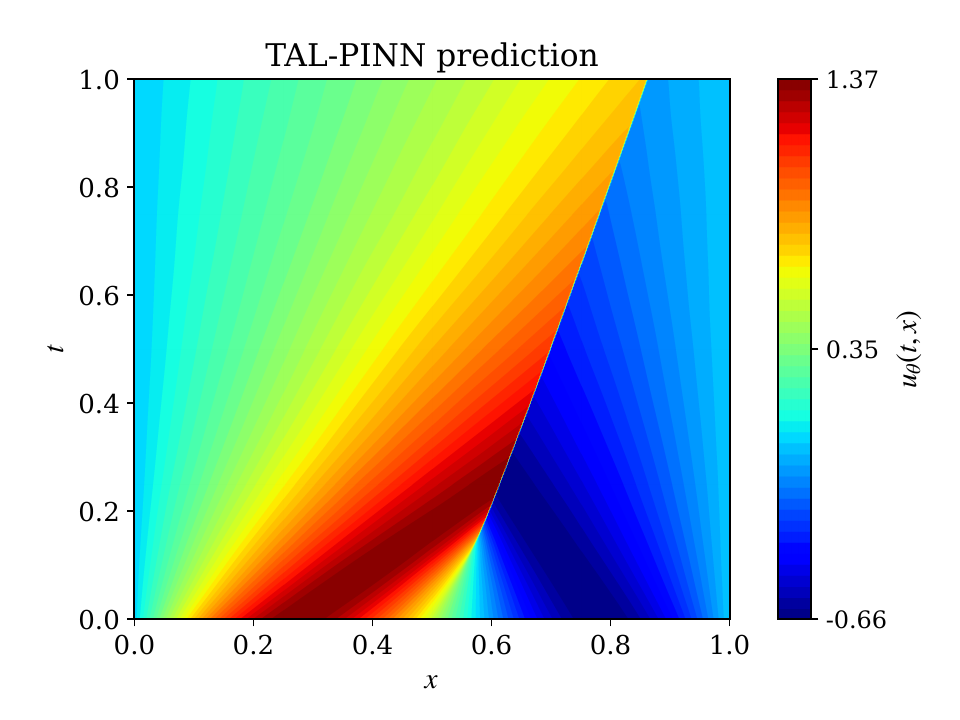}
    \caption{TAL-PINN prediction}
\end{subfigure}\hfill
\begin{subfigure}[t]{0.480\linewidth}
    \centering
    \includegraphics[width=\linewidth]{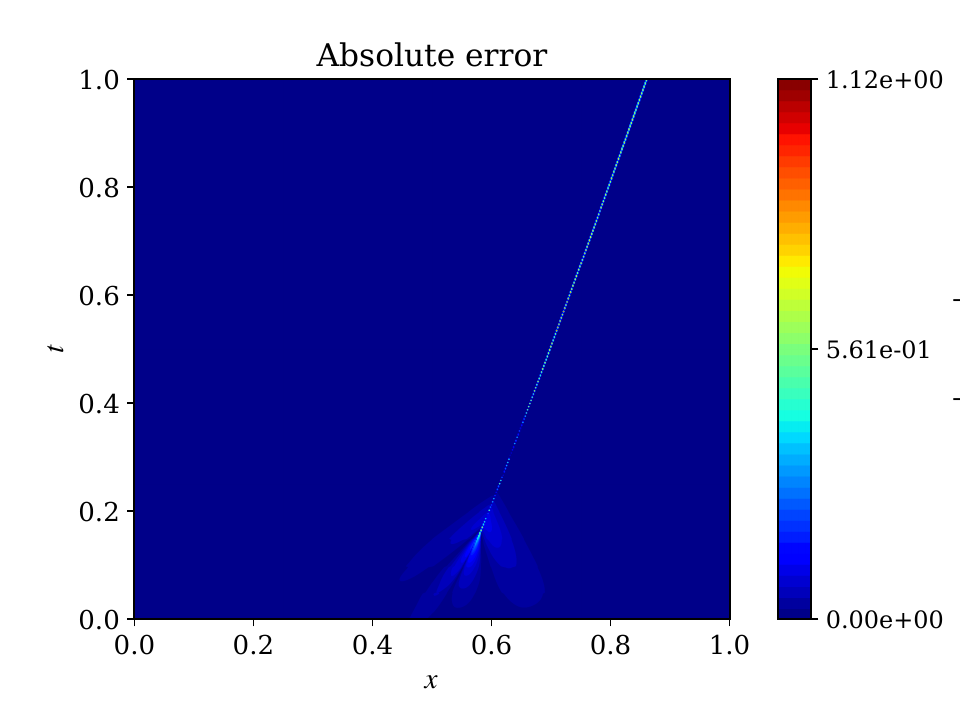}
    \caption{Absolute error}
\end{subfigure}

\par\medskip

\begin{subfigure}[t]{0.480\linewidth}
    \centering
    \includegraphics[width=\linewidth]{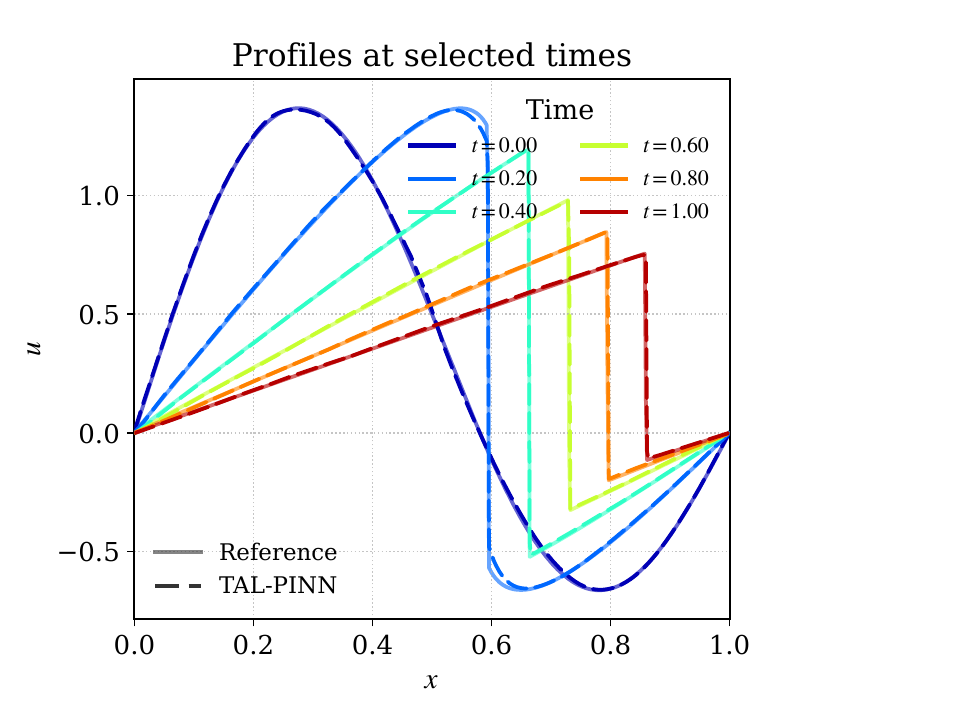}
    \caption{Time profiles}
\end{subfigure}\hfill
\begin{subfigure}[t]{0.480\linewidth}
    \centering
    \includegraphics[width=\linewidth]{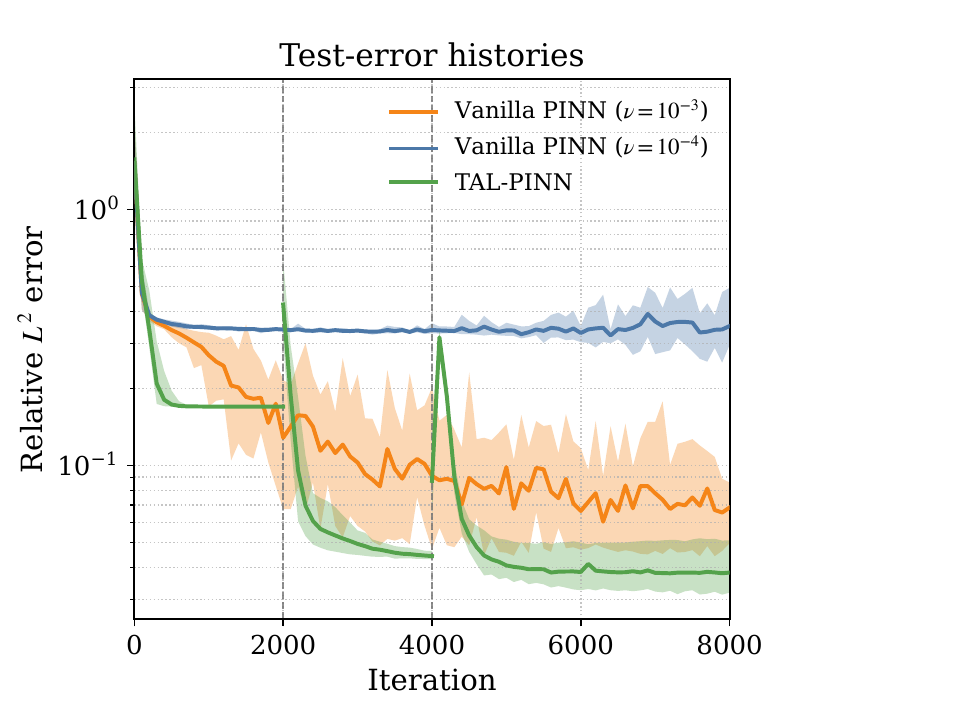}
    \caption{Test-error history}
\end{subfigure}
\caption{Moving Burgers shock: (a) inviscid reference solution; (b) adaptive mesh; (c) TAL-PINN prediction; (d) absolute error; (e) profiles at selected times; and (f) relative $L^2$ test-error histories. Panels (b)--(e) show results from one run. Panel (f) shows mean errors with minimum--maximum shading; TAL-PINN changes viscosity after 2000 and 4000 PDE iterations.}
\label{fig:exp_moving_solution}
\end{figure}

Figure~\ref{fig:exp_moving_solution}(a,c,d) presents the reference solution, TAL-PINN prediction, and absolute error. TAL-PINN captures the curved shock trajectory together with the surrounding smooth solution.
Table~\ref{tab:exp_moving} reports a final relative $L^2$ error of $3.82\times10^{-2}$, compared with $3.50\times10^{-1}$ for vanilla PINN at the same final viscosity $10^{-4}$. TAL-PINN also improves upon the larger-viscosity baseline.

\begin{table}[!htb]
\centering
\caption{Moving Burgers shock: mean final relative $L^2$ errors.}
\label{tab:exp_moving}
\small
\begin{tabular}{@{}lcc@{}}
\toprule
Method & Final $\nu$ & Relative $L^2$\\
\midrule
Vanilla PINN & $10^{-3}$ & $6.85\times10^{-2}$\\
Vanilla PINN & $10^{-4}$ & $3.50\times10^{-1}$\\
TAL-PINN & $10^{-4}$ & $3.82\times10^{-2}$\\
\bottomrule
\end{tabular}
\end{table}

The mesh in Figure~\ref{fig:exp_moving_solution}(b) follows the shock trajectory, while the predicted profiles in panel (e) reproduce its location and steep transition at successive times. The error histories in panel (f) show that TAL-PINN achieves further accuracy gains through the successive viscosity stages, whereas direct training at the smallest viscosity levels off at a substantially larger error. These results demonstrate that the adaptive lifting strategy remains effective for moving shocks without prescribing the shock trajectory.

\subsection{1D Euler equations}
\label{sec:euler_1d}

We next consider the one-dimensional compressible Euler equations,
\begin{equation}
    \mathbf{Q}_t+\mathbf{F}(\mathbf{Q})_x=0,\qquad
    \mathbf{Q}=\begin{pmatrix}\rho\\\rho u\\E\end{pmatrix},\qquad
    \mathbf{F}=\begin{pmatrix}\rho u\\\rho u^2+p\\u(E+p)\end{pmatrix},
    \label{eq:exp_euler}
\end{equation}
where $E=p/(\gamma-1)+\rho u^2/2$ and $\gamma=1.4$. The physical domain is $x\in[0,1]$, with a Riemann initial condition separated at $x=0.5$. In training, a constant artificial viscosity is applied to the conservative variables,
\[
    \mathbf{Q}_t+\mathbf{F}(\mathbf{Q})_x=\nu\mathbf{Q}_{xx}.
\]
The loss combines the three conservative residuals and the initial conditions. No additional boundary loss is used over the time intervals considered here. The inviscid reference is obtained from the exact Riemann solution~\cite{toro2013riemann}. Relative $L^2$ errors are evaluated componentwise over the full $129\times257$ space--time grid, while the profiles below are taken at the final time.

The solution network has a shared Fourier feature layer and a two-layer trunk, followed by three independent branches for $(\rho,u,p)$. The first trunk layer and the two hidden layers of each branch use SiLU, and the second trunk layer is linear. Each branch has a scalar output without bias. Positivity of $\rho$ and $p$ is enforced using softplus.

For both shock tubes, the monitor emphasizes compression in the velocity field. With $g=|u_x|-u_x$, we use
\[
    M(t,x)=1+\frac{g(t,x)}{2\alpha(t)},\qquad
    \alpha(t)=\max\!\left(\int_0^1 g(t,x)\,\dx,\,1\right).
\]
After smoothing the monitor in space, we construct the adaptive coordinate using \eqref{eq:equidistribution_1d_map} at each time. The coordinate loss uses $w_J=0.1$ and $\delta_J=0.05$. One mesh update and one coordinate fit separate the two PDE stages.

\subsubsection{Sod shock tube}
\label{sec:exp_sod}

For the Sod problem, the initial states are
\[
(\rho,u,p)(0,x)=
\begin{cases}
(1,0,1),&x<0.5,\\
(0.125,0,0.1),&x\geq0.5,
\end{cases}
\]
and the final time is $T_{\mathrm{end}}=0.2$. TAL-PINN reduces the viscosity from $10^{-3}$ to $10^{-4}$, while the vanilla PINN uses $\nu=10^{-4}$ throughout.

\begin{figure}[!htbp]
\centering

\begin{subfigure}[t]{0.480\linewidth}
    \centering
    \includegraphics[width=\linewidth]{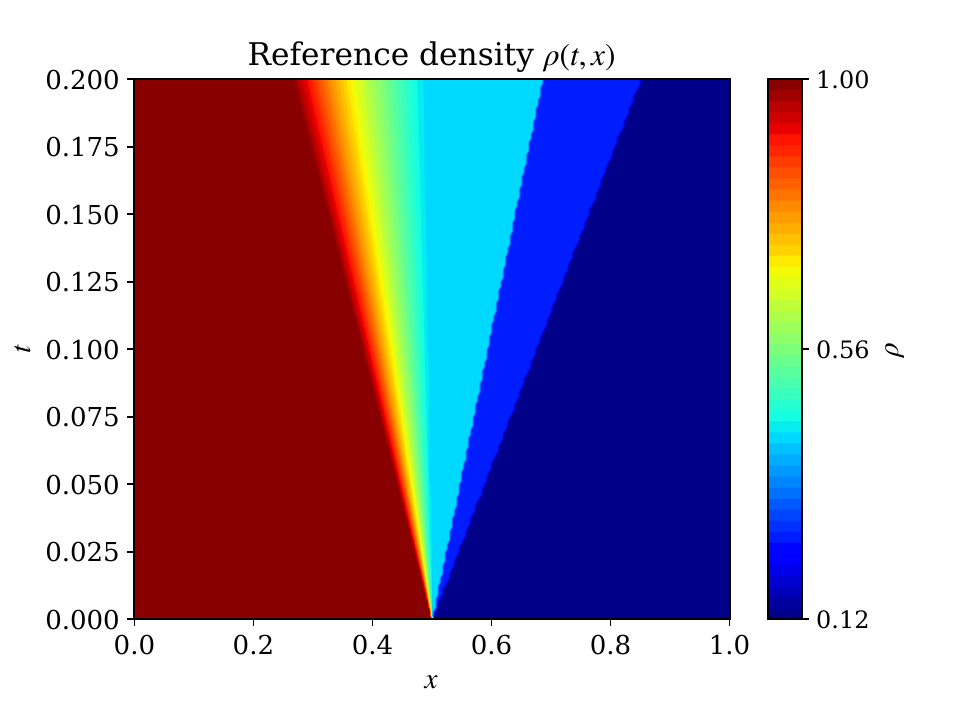}
    \caption{Reference density}
\end{subfigure}\hfill
\begin{subfigure}[t]{0.480\linewidth}
    \centering
    \includegraphics[width=\linewidth]{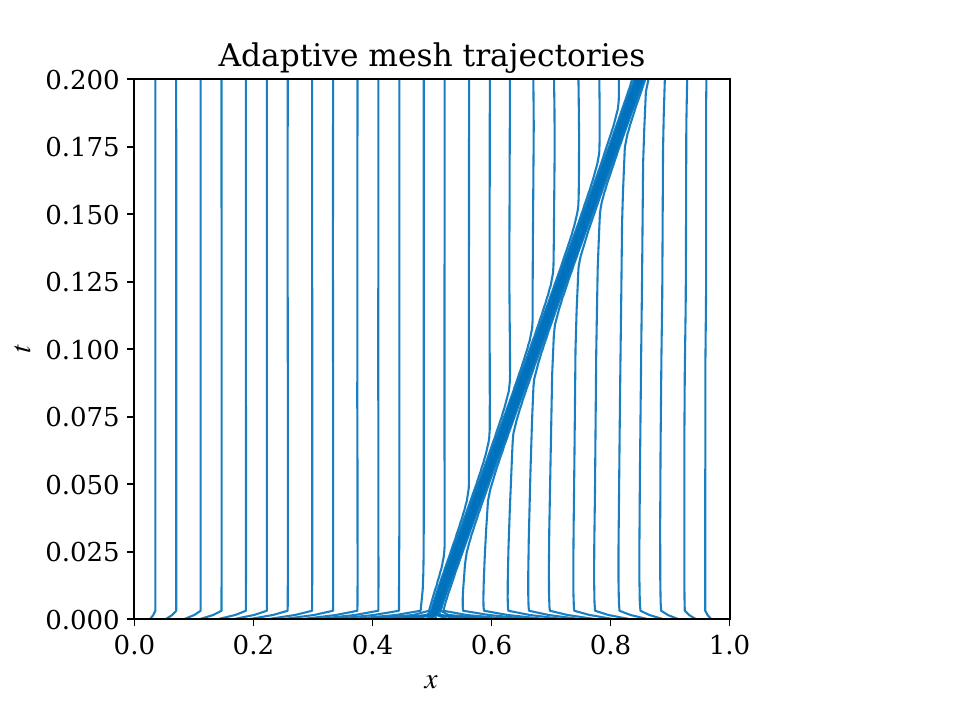}
    \caption{Adaptive mesh}
\end{subfigure}

\par\medskip

\begin{subfigure}[t]{0.480\linewidth}
    \centering
    \includegraphics[width=\linewidth]{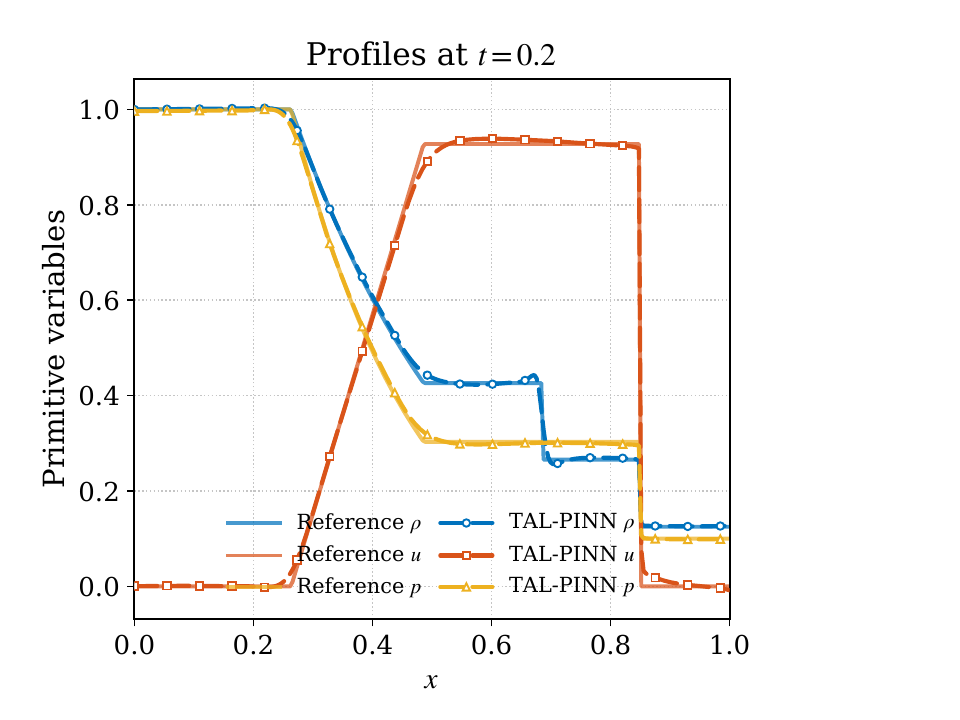}
    \caption{Profiles at $t=0.2$}
\end{subfigure}\hfill
\begin{subfigure}[t]{0.480\linewidth}
    \centering
    \includegraphics[width=\linewidth]{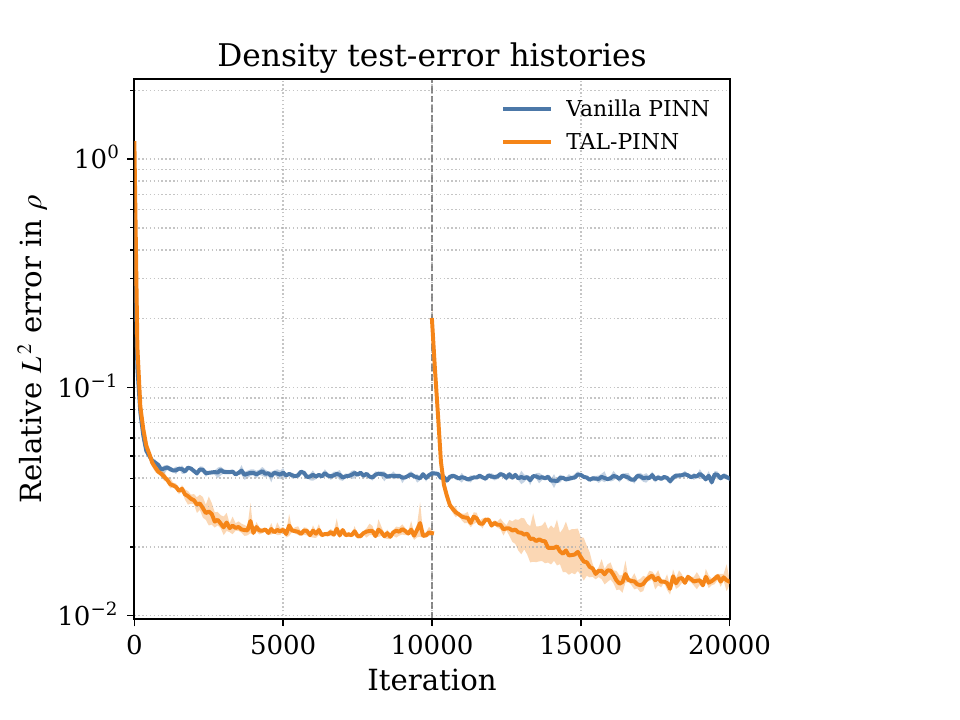}
    \caption{Density error}
\end{subfigure}

\par\medskip

\begin{subfigure}[t]{0.480\linewidth}
    \centering
    \includegraphics[width=\linewidth]{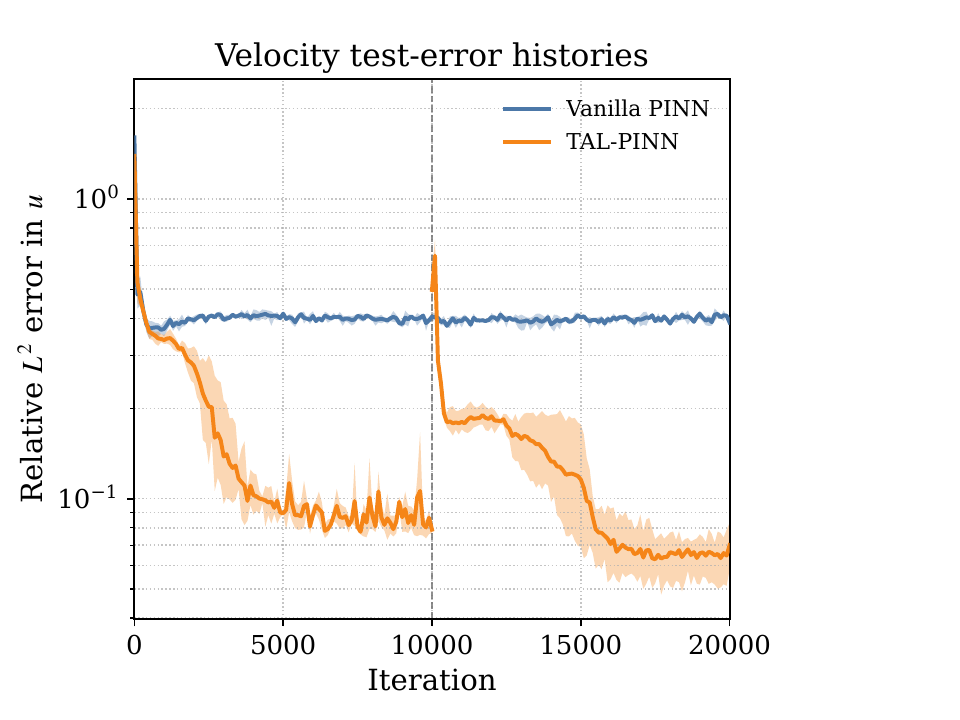}
    \caption{Velocity error}
\end{subfigure}\hfill
\begin{subfigure}[t]{0.480\linewidth}
    \centering
    \includegraphics[width=\linewidth]{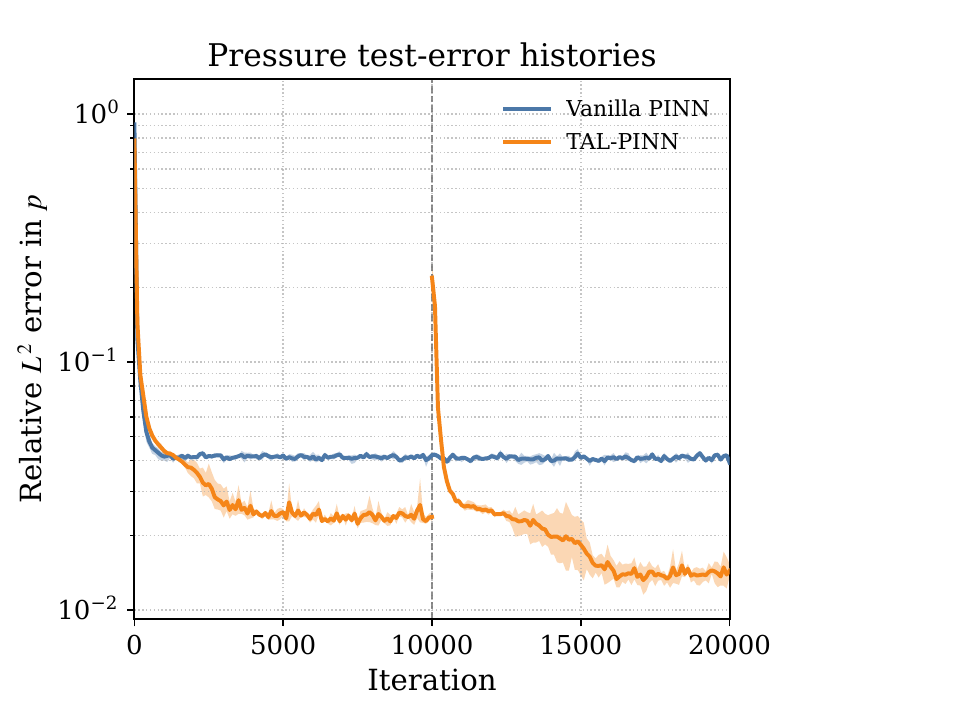}
    \caption{Pressure error}
\end{subfigure}

\caption{Sod shock tube: (a) reference density; (b) adaptive mesh; (c) profiles at $t=0.2$; and (d)--(f) relative $L^2$ test-error histories for density, velocity, and pressure, respectively. Panels (b,c) show results from one run. In (c), solid curves denote the inviscid reference and markers denote TAL-PINN predictions; colors distinguish the three variables. Panels (d)--(f) show mean errors with minimum--maximum shading. TAL-PINN introduces the fitted coordinate and reduces the viscosity after 10000 PDE iterations.}
\label{fig:exp_sod_convergence}
\end{figure}

Figure~\ref{fig:exp_sod_convergence}(a,b) shows that the adaptive mesh constructed from the coarse-stage velocity monitor concentrates along the right-moving shock trajectory. With the fitted coordinate as the lifting input, TAL-PINN accurately reproduces the Riemann wave pattern in density, velocity, and pressure, including the steep front and the intermediate states, as shown by the final-time profiles in panel (c).

The componentwise error histories in panels (d)--(f) show further error reduction in the second stage, after introducing the fitted coordinate and lowering the viscosity. All three variables reach smaller errors than the vanilla baseline. As reported in Table~\ref{tab:exp_euler_errors}, the velocity relative $L^2$ error is $7.00\times10^{-2}$ for TAL-PINN, compared with $3.87\times10^{-1}$ for the vanilla PINN. These results demonstrate that TAL-PINN improves the approximation of the coupled system with the same total number of PDE updates as the vanilla baseline.

\subsubsection{Lax shock tube}
\label{sec:euler_1d_lax}

For the Lax problem, the initial states are
\[
(\rho,u,p)(0,x)=
\begin{cases}
(0.445,0.698,3.528),&x<0.5,\\
(0.5,0,0.571),&x\geq0.5,
\end{cases}
\]
with $T_{\mathrm{end}}=0.14$. TAL-PINN reduces the viscosity from $3\times10^{-3}$ to $3\times10^{-4}$, while the vanilla baseline uses $\nu=3\times10^{-4}$ throughout.

\begin{figure}[!htbp]
\centering

\begin{subfigure}[t]{0.480\linewidth}
    \centering
    \includegraphics[width=\linewidth]{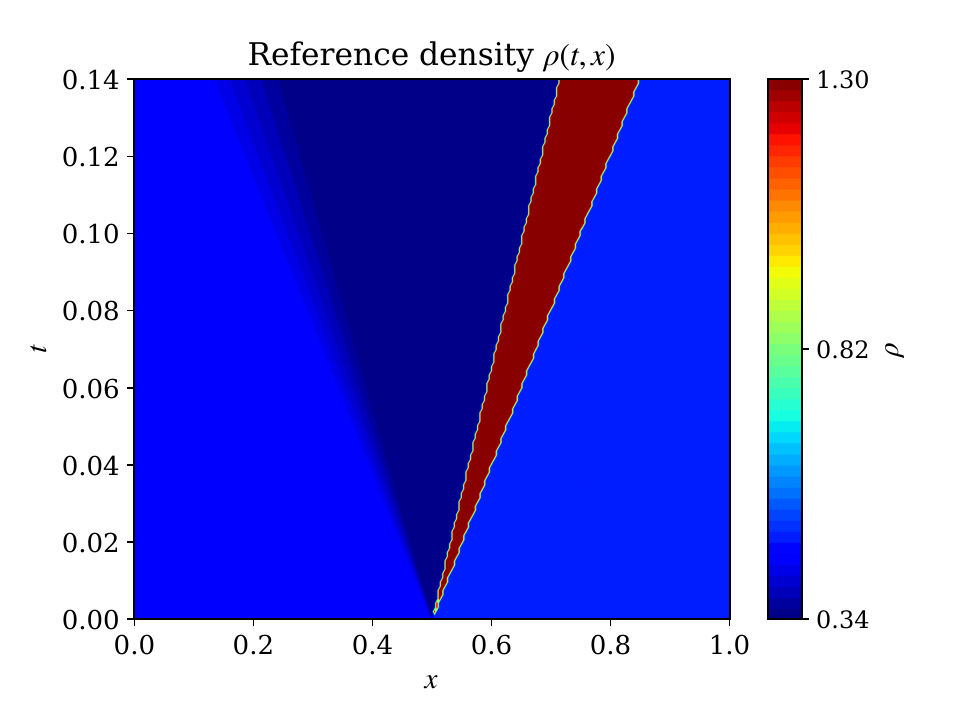}
    \caption{Reference density}
\end{subfigure}\hfill
\begin{subfigure}[t]{0.480\linewidth}
    \centering
    \includegraphics[width=\linewidth]{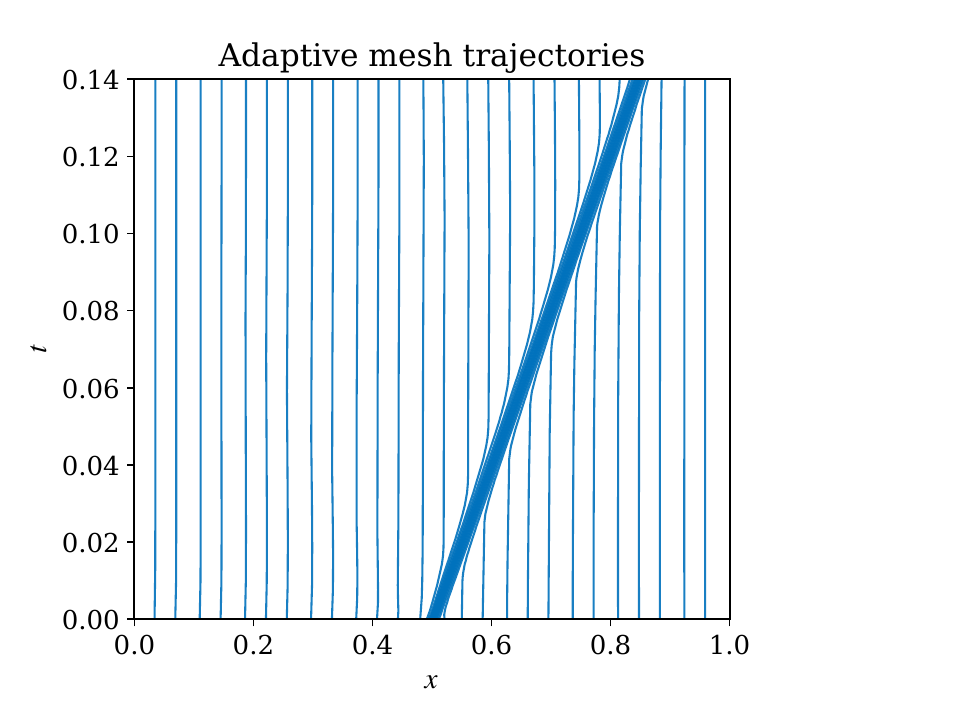}
    \caption{Adaptive mesh}
\end{subfigure}

\par\medskip

\begin{subfigure}[t]{0.480\linewidth}
    \centering
    \includegraphics[width=\linewidth]{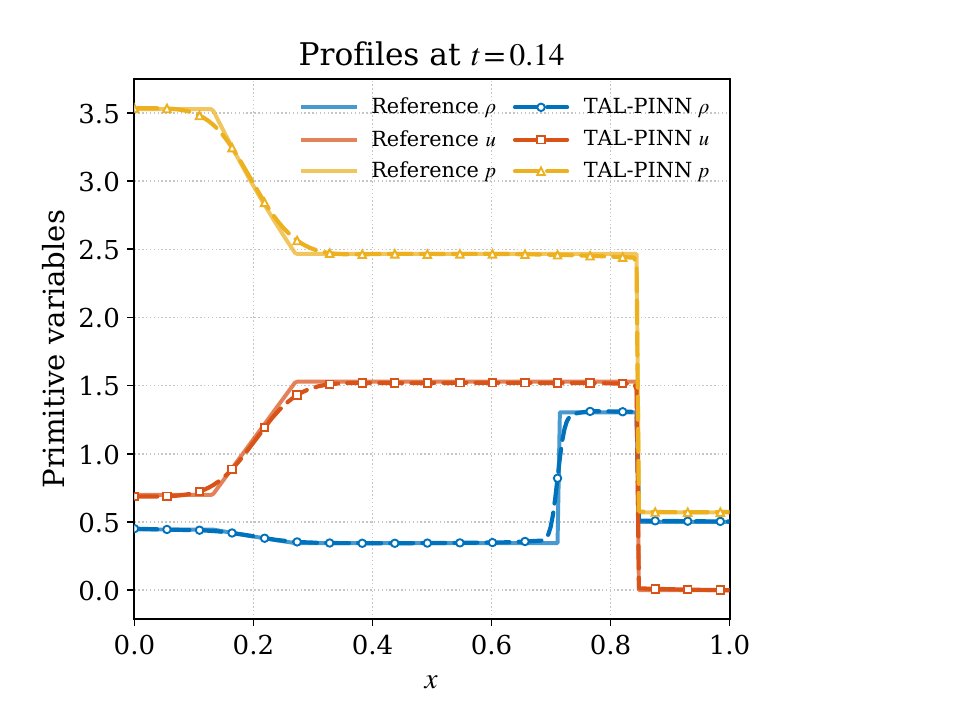}
    \caption{Profiles at $t=0.14$}
\end{subfigure}\hfill
\begin{subfigure}[t]{0.480\linewidth}
    \centering
    \includegraphics[width=\linewidth]{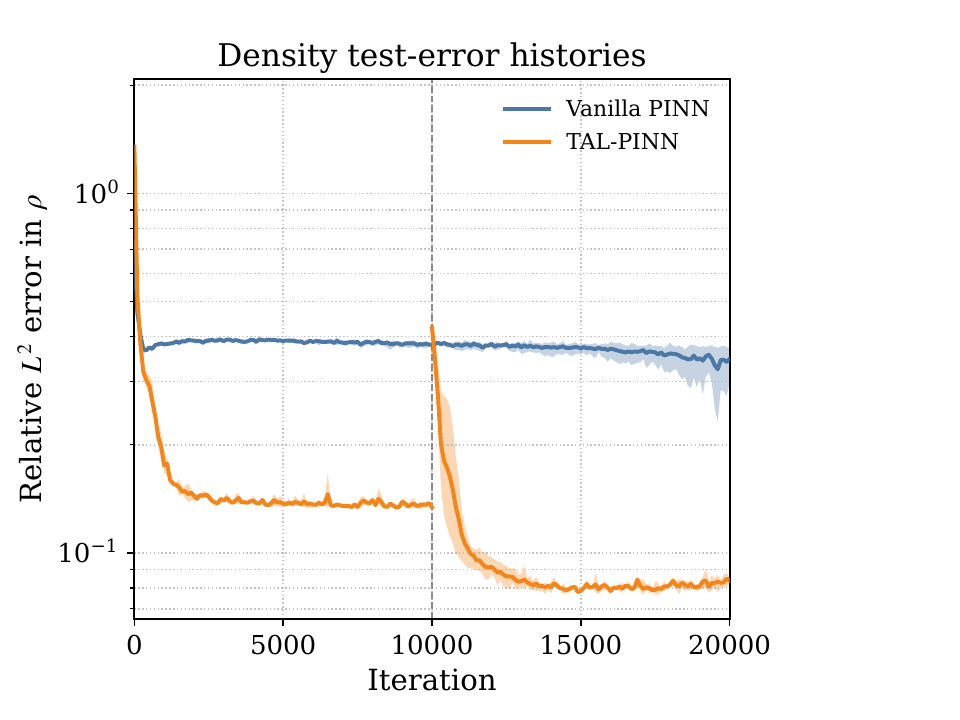}
    \caption{Density error}
\end{subfigure}

\par\medskip

\begin{subfigure}[t]{0.480\linewidth}
    \centering
    \includegraphics[width=\linewidth]{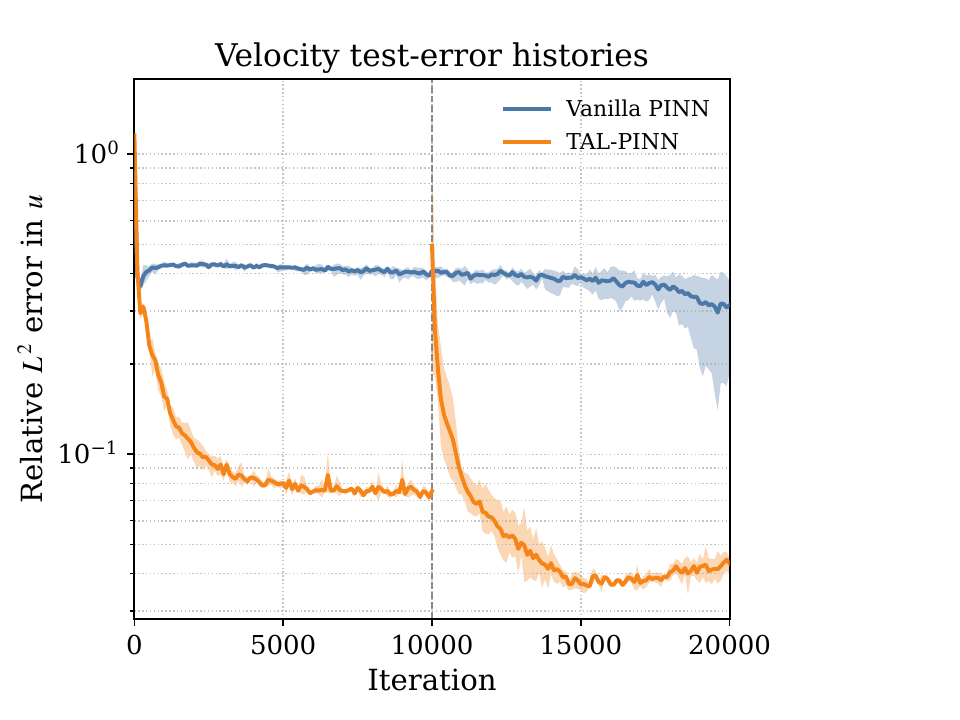}
    \caption{Velocity error}
\end{subfigure}\hfill
\begin{subfigure}[t]{0.480\linewidth}
    \centering
    \includegraphics[width=\linewidth]{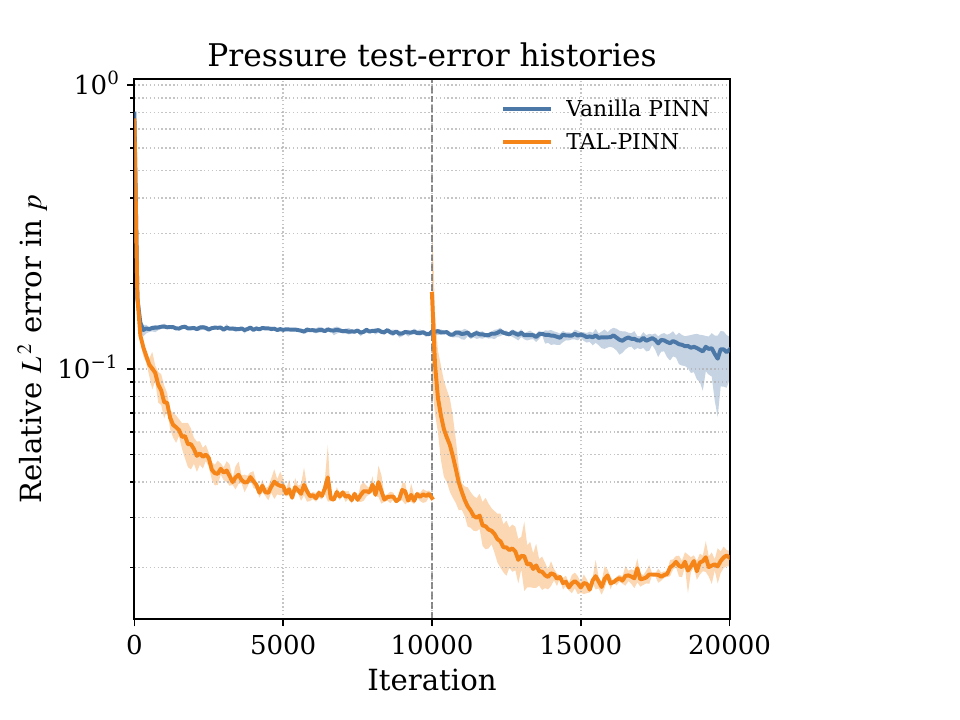}
    \caption{Pressure error}
\end{subfigure}

\caption{Lax shock tube: (a) reference density; (b) adaptive mesh; (c) profiles at $t=0.14$; and (d)--(f) relative $L^2$ test-error histories for density, velocity, and pressure, respectively. Panels (b,c) show results from one run. In (c), solid curves denote the inviscid reference and markers denote TAL-PINN predictions; colors distinguish the three variables. Panels (d)--(f) show mean errors with minimum--maximum shading. TAL-PINN introduces the fitted coordinate and reduces the viscosity after 10000 PDE iterations.}
\label{fig:exp_lax_convergence}
\end{figure}

The Lax problem produces a pronounced density peak and substantial variations in all three solution components. Figure~\ref{fig:exp_lax_convergence}(a,b) shows that the adaptive mesh follows the moving shock. The final-time profiles in panel (c) demonstrate that TAL-PINN resolves the steep wave structure and reproduces the density peak together with the intermediate velocity and pressure plateaus.

The componentwise error histories in panels (d)--(f) show a pronounced error reduction in the second stage, while the vanilla PINN remains at higher error levels. Using the coordinate constructed from the velocity monitor, TAL-PINN improves all three solution components relative to the vanilla baseline, as reported in Table~\ref{tab:exp_euler_errors}. The Sod and Lax results therefore extend the accuracy and training benefits observed for scalar shocks to a coupled hyperbolic system.

\begin{table}[!htb]
\centering
\caption{Euler shock tubes: mean componentwise relative $L^2$ errors on the full space--time test grid. The two methods use the same final viscosity within each problem: $10^{-4}$ for Sod and $3\times10^{-4}$ for Lax.}
\label{tab:exp_euler_errors}
\small
\setlength{\tabcolsep}{4pt}
\begin{tabular}{@{}llccc@{}}
\toprule
Problem & Method & $\rho$ & $u$ & $p$\\
\midrule
Sod & Vanilla PINN & $3.98\times10^{-2}$ & $3.87\times10^{-1}$ & $3.93\times10^{-2}$\\
 & TAL-PINN & $1.41\times10^{-2}$ & $7.00\times10^{-2}$ & $1.44\times10^{-2}$\\
\midrule
Lax & Vanilla PINN & $3.46\times10^{-1}$ & $3.14\times10^{-1}$ & $1.17\times10^{-1}$\\
 & TAL-PINN & $8.39\times10^{-2}$ & $4.33\times10^{-2}$ & $2.17\times10^{-2}$\\
\bottomrule
\end{tabular}
\end{table}

\subsection{2D Burgers equation}
\label{sec:burgers_2d}

To examine the extension to two spatial dimensions, we consider
\begin{equation}
    u_t+\left(\frac{u^2}{2}\right)_x+\left(\frac{u^2}{2}\right)_y=0,
    \qquad (t,x,y)\in[0,1]\times[0,4]^2,
    \label{eq:exp_burgers2d}
\end{equation}
with periodic boundary conditions and $u(0,x,y)=\sin(\pi(x+y)/2)$. The inviscid solution develops oblique shock lines. The training equation uses the Laplacian regularization
\[
    u_t+u(u_x+u_y)=\nu(u_{xx}+u_{yy}).
\]
The reference uses a fifth-order finite-difference WENO-Z scheme~\cite{borges2008wenoz} and third-order SSP Runge--Kutta time integration. The relative $L^2$ error is evaluated at $t=1$ on a $257\times257$ spatial grid.

TAL-PINN uses two viscosity stages, $10^{-2}$ and $10^{-3}$, while the vanilla model uses $\nu=10^{-3}$ throughout. The lifted input is $(t,x,y,\xi,\eta)$, with $(\xi,\eta)$ learned from $(t,x,y)$. Opposite-edge value matching enforces periodicity. With $D=u_x+u_y$ and $g=|D|-D$, the monitor is
\[
    M=1+\frac{g}{2\alpha(t)},\qquad
    \alpha(t)=\max\!\left(\frac{1}{|\Omega_P|}\int_{\Omega_P}g\,\dx\dy,\,1\right).
\]
The adaptive mesh is generated by an MMPDE5-type evolution initialized from the elliptic mesh condition \eqref{eq:winslow}. The coordinate loss uses $w_J=0.1$ and $\delta_J=10^{-3}$.

\begin{figure}[!htb]
\centering
\begin{subfigure}[t]{0.480\linewidth}
    \centering
    \includegraphics[width=\linewidth]{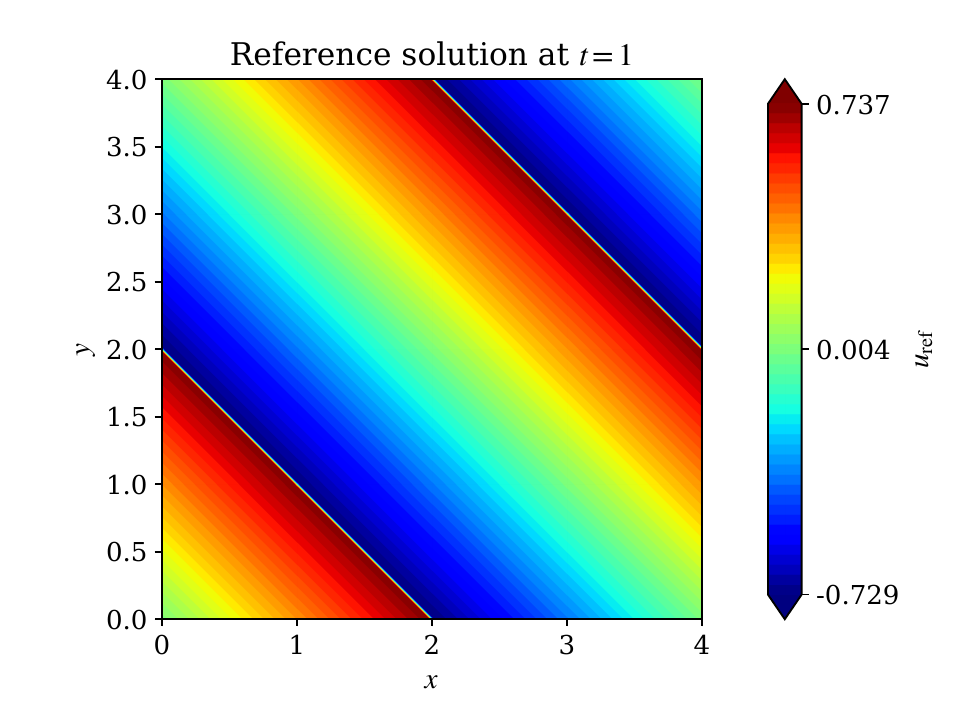}
    \caption{Inviscid reference}
\end{subfigure}
\hfill
\begin{subfigure}[t]{0.480\linewidth}
    \centering
    \includegraphics[width=\linewidth]{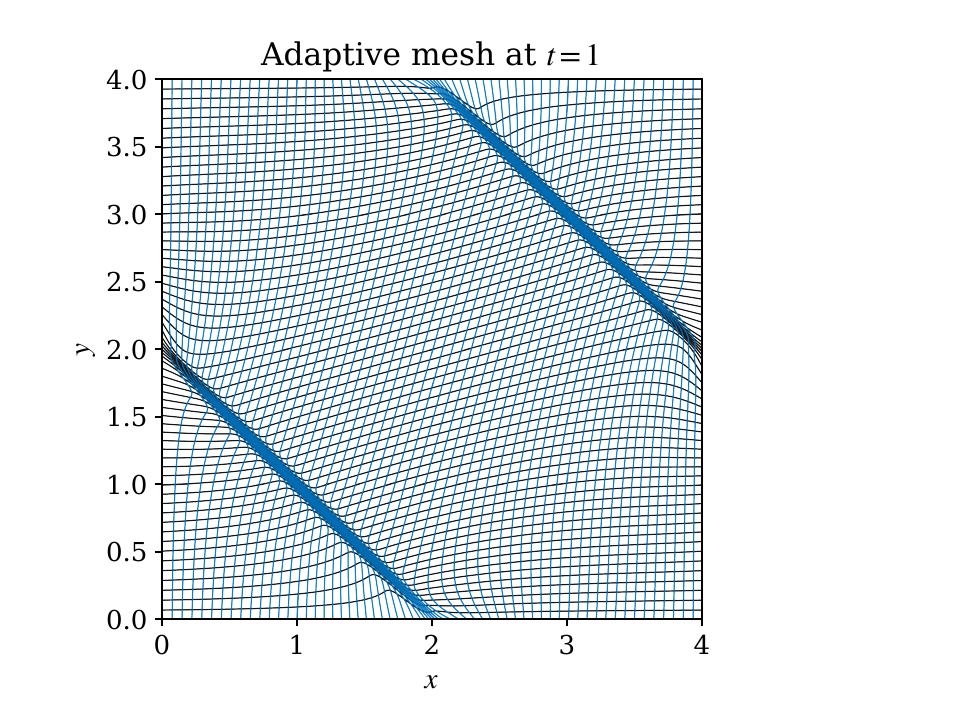}
    \caption{Adaptive mesh}
\end{subfigure}
\par\medskip
\begin{subfigure}[t]{0.480\linewidth}
    \centering
    \includegraphics[width=\linewidth]{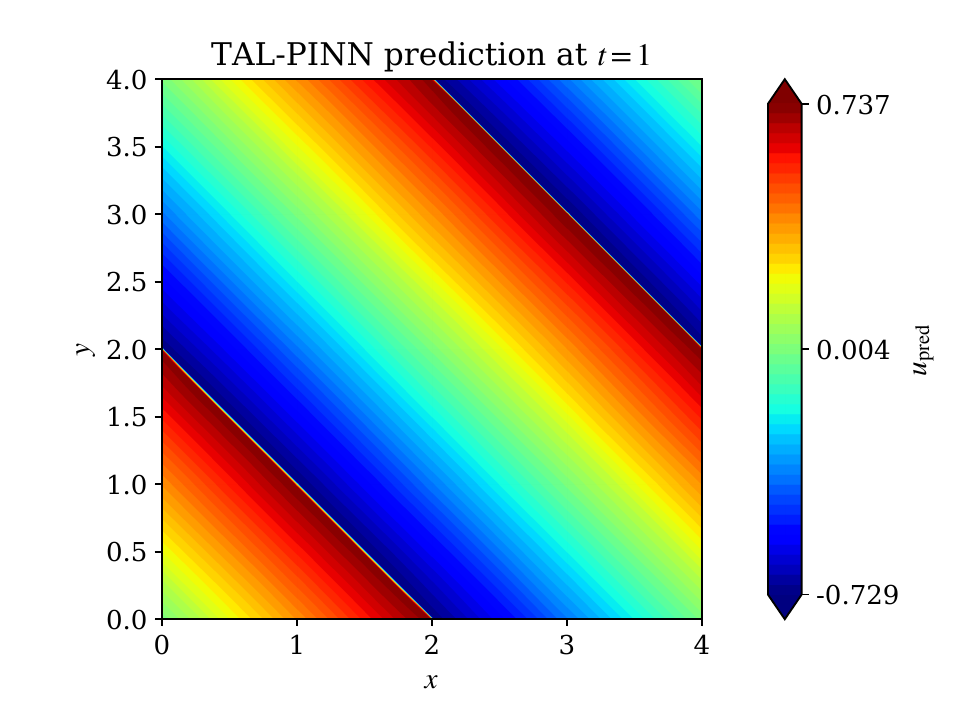}
    \caption{TAL-PINN prediction}
\end{subfigure}
\hfill
\begin{subfigure}[t]{0.480\linewidth}
    \centering
    \includegraphics[width=\linewidth]{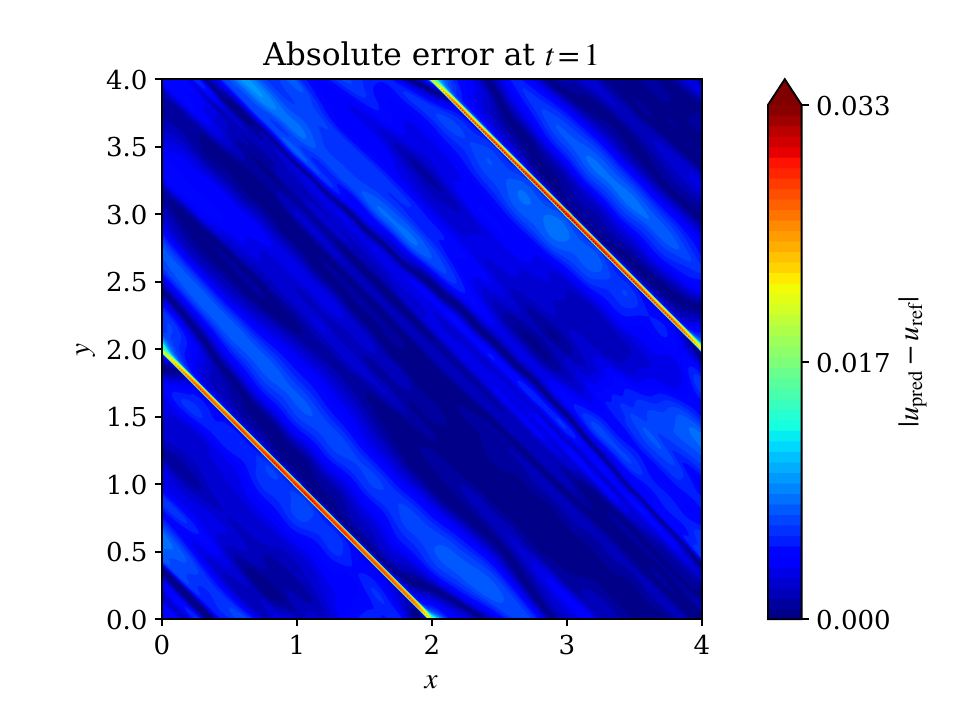}
    \caption{Absolute error}
\end{subfigure}
\caption{2D Burgers equation at $t=1$: reference solution, adaptive mesh, TAL-PINN prediction, and absolute error from one run. The reference and prediction share a color scale. Panel (b) displays every fourth mesh line.}
\label{fig:exp_burgers2d_solution}
\end{figure}

Figure~\ref{fig:exp_burgers2d_solution}(a,c,d) presents the reference solution, TAL-PINN prediction, and absolute error at the final time. TAL-PINN captures the oblique shock lines and preserves the smooth variation between them, with the remaining error concentrated near the steep fronts. Table~\ref{tab:exp_burgers2d} reports a relative $L^2$ error of $1.47\times10^{-2}$, compared with $4.55\times10^{-1}$ for vanilla PINN using the same 14000 PDE iterations.

\begin{table}[!htb]
\centering
\caption{2D Burgers equation: mean final-time relative $L^2$ errors at $\nu=10^{-3}$.}
\label{tab:exp_burgers2d}
\small
\begin{tabular}{@{}lc@{}}
\toprule
Method & Relative $L^2$\\
\midrule
Vanilla PINN & $4.55\times10^{-1}$\\
TAL-PINN & $1.47\times10^{-2}$\\
\bottomrule
\end{tabular}
\end{table}

\begin{figure}[!htb]
\centering
\includegraphics[width=0.48\linewidth]{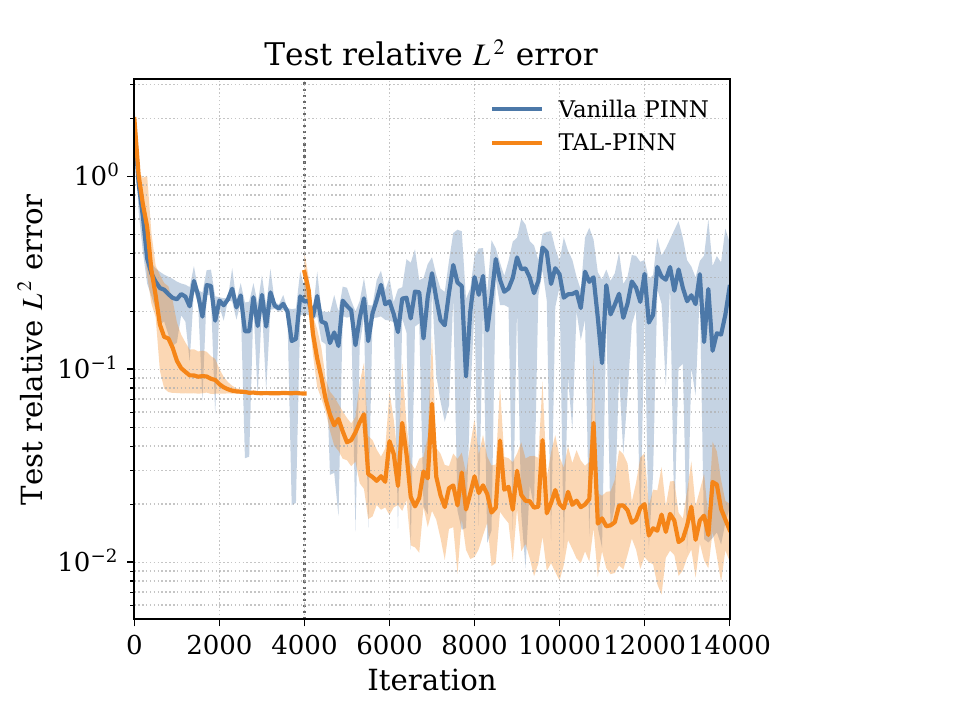}
\caption{2D Burgers equation: final-time relative $L^2$ test-error histories. Shading indicates the minimum--maximum range; the coordinate is updated and the viscosity is reduced after 4000 PDE iterations.}
\label{fig:exp_burgers2d_convergence}
\end{figure}

The adaptive mesh in Figure~\ref{fig:exp_burgers2d_solution}(b) concentrates resolution along the oblique shocks, while Figure~\ref{fig:exp_burgers2d_convergence} shows a substantial accuracy gain after the coordinate update. The separation from the vanilla baseline persists across repeated runs. This example demonstrates that the same lifting principle remains effective with two spatial coordinates and a larger residual sampling budget, enabling accurate small-viscosity training for a multidimensional shock pattern.

\subsection{Computational cost}
\label{sec:exp_cost}

We next examine the computational cost associated with the accuracy improvements reported above. Table~\ref{tab:exp_cost} reports the time spent on solution-network training, adaptive mesh generation, and coordinate-network fitting.

\begin{table}[!htb]
\centering
\caption{Mean computational costs in seconds, accumulated across all stages and updates. A dash denotes an operation not used by the method.}
\label{tab:exp_cost}
\small
\setlength{\tabcolsep}{3pt}
\begin{tabular*}{\linewidth}{@{\extracolsep{\fill}}llrrrr@{}}
\toprule
Problem & Method (final $\nu$) & PDE & Mesh & Coordinate & Total\\
\midrule
Stationary & Vanilla ($10^{-3}$) & 24.29 & --- & --- & 24.29\\
Burgers & Vanilla ($10^{-4}$) & 23.97 & --- & --- & 23.97\\
 & TAL-PINN ($10^{-4}$) & 28.89 & 0.93 & 38.74 & 68.56\\
 & VC ($10^{-4}$) & 36.68 & --- & --- & 36.68\\
 & AS+VC ($10^{-4}$) & 41.10 & 1.07 & --- & 42.17\\
\midrule
Moving & Vanilla ($10^{-3}$) & 23.85 & --- & --- & 23.85\\
Burgers & Vanilla ($10^{-4}$) & 24.16 & --- & --- & 24.16\\
 & TAL-PINN ($10^{-4}$) & 28.93 & 0.96 & 39.50 & 69.39\\
\midrule
Sod & Vanilla ($10^{-4}$) & 198.22 & --- & --- & 198.22\\
 & TAL-PINN ($10^{-4}$) & 255.16 & 1.17 & 20.29 & 276.62\\
Lax & Vanilla ($3\times10^{-4}$) & 196.62 & --- & --- & 196.62\\
 & TAL-PINN ($3\times10^{-4}$) & 254.92 & 1.16 & 20.14 & 276.22\\
\midrule
2D Burgers & Vanilla ($10^{-3}$) & 126.54 & --- & --- & 126.54\\
 & TAL-PINN ($10^{-3}$) & 198.55 & 22.21 & 24.48 & 245.24\\
\bottomrule
\end{tabular*}
\end{table}

Relative to vanilla PINNs at the same final viscosity, TAL-PINN requires approximately 2.9 times the total runtime for the one-dimensional Burgers tests, 1.4 times for Euler, and 1.9 times for two-dimensional Burgers. These increases accompany the substantial accuracy gains reported above. For the stationary shock, the longer PDE-training schedules of VC and AS+VC still do not recover TAL-PINN's small-viscosity accuracy.

Coordinate fitting accounts for a substantial fraction of the total cost in the one-dimensional Burgers tests, where two fits are performed and solution-network training is relatively inexpensive. With the same fitting batch size and iteration budget, the cost per fit increases only modestly from one to two spatial dimensions, much less than the increase in solution-network training time. Coordinate fitting is a supervised regression task on known nodal data, whereas solution-network training requires evaluating PDE residuals involving space--time derivatives and second-order viscous terms. The larger residual batches and longer PDE-training schedule in the two-dimensional test further increase the solution-training cost.

Mesh generation contributes little to the total cost in the one-dimensional tests and becomes more expensive in two dimensions, although its cost remains well below the solution-network training time. In the Euler and two-dimensional Burgers tests, solution-network training accounts for most of the total cost. The adaptive mesh and fitted coordinate network are constructed between stages and used throughout the subsequent solution-network training.

\section{Conclusion \& Discussion}
\label{sec:conclusion}

TAL-PINN constructs solution-dependent lifting coordinates from a coarse viscous approximation and uses them to train at smaller viscosity. The adaptive coordinates encode the transition geometry in smooth auxiliary inputs and concentrate physical samples near steep layers, with measure correction preserving the physical residual objective. This construction integrates naturally with viscosity continuation and replaces manually prescribed interface indicators. Our analysis establishes a scalar a posteriori $L^2$ error decomposition into viscous approximation error, statistical error, and empirical training loss, together with a variance-reduction criterion for adaptive residual sampling. The NTK analysis relates the additional coordinate variation to enriched tangent features and faster residual-mode decay.

Numerical experiments on benchmark problems---including stationary and moving 1D Burgers shocks, the 1D Euler Sod and Lax shock tubes, and a 2D Burgers problem---demonstrate that TAL-PINN resolves steep shock structures and improves accuracy with the same PDE-training budget as vanilla PINNs. The stationary-shock ablations show that lifting enables continued accuracy gains as viscosity decreases, beyond those obtained by viscosity continuation and adaptive sampling alone. The timing results show that coordinate fitting accounts for a substantial share of the cost in the one-dimensional Burgers tests, while solution-network training dominates in the Euler and two-dimensional problems.

Extending this autonomous lifting mechanism to more complex geometries calls for efficient construction of the full-domain $r$-adaptive mapping. A natural direction for future work is to integrate the TAL-PINN strategy with domain decomposition (e.g., FBPINNs~\cite{moseley2023fbpinns} or XPINNs~\cite{jagtap2020xpinns}). In such a coupled framework, the data-driven monitor function could activate lifting within subdomains containing steep gradients, with the aim of reducing geometric preparation costs in large-scale three-dimensional applications.

\section*{Acknowledgement}

We are grateful to the referees for their insightful and constructive comments, which have significantly improved the paper.
This work was supported by the National Key R\&D Program of China (2024YFA1012600), and by the NSF of China grant 12171237.

\section*{Code availability}

The implementation of TAL-PINN and the code for the numerical experiments are available at \url{https://github.com/wszhuym/TAL-PINN}.



\end{document}